\documentclass{article}
\usepackage{amssymb,amsmath,titlesec}
\usepackage[all]{xy}

\titleformat{\section}[block]
  {\sc\large\filcenter}
  {\S\thesection.}{.5em}{}
\titleformat{\subsection}[runin]
  {\normalfont\bfseries}
  {\thesubsection.}{.5em}{}

\def\t{\textrm}
\def\fn{\footnotesize}
\def\bb{\mathbb}
\def\mc{\mathcal}
\def\mf{\mathfrak}
\def\beq{\begin{eqnarray}}
\def\eeq{\end{eqnarray}}
\def\beqq{\begin{eqnarray*}}
\def\eeqq{\end{eqnarray*}}

\def\CDO{\mc{D}^{\t{ch}}}
\def\CDC{\E^{\t{ch},*}}
\def\CDCnew{\wt\E^{\t{ch},*}}
\def\d{\partial}
\def\E{\mc{E}}
\def\A{\hat{A}}
\def\dbar{\bar{\partial}}
\def\dbardef{\left[\begin{array}{cc} \dbar & \bar\Delta \\ 0 & \bar\d
  \end{array}\right]}
\def\Tr{\t{Tr}\;}
\def\T{\mc{T}}
\def\wh{\widehat}
\def\wt{\widetilde}
\def\otimesez{\otimes_{\t{ez}}}

\title{\sf THE WITTEN GENUS AND VERTEX ALGEBRAS}
\author{\sf Pokman Cheung}
\date{\sf February 15, 2010}

\begin{document}

\maketitle

\begin{abstract}
This article is the first report of an ongoing project aimed at finding 
a geometric interpretation of the Witten genus and other tmf classes.
Section 2 reviews the sheaves of chiral differential operators 
$\CDO_{M,\xi}$ over a complex manifold $M$, including their construction, 
obstructions and relation with the Witten genus of $M$.
In section 3, the structure of $\CDO_{M,\xi}$ as a sheaf of vertex 
algebras is reorganized in terms of $\mc{O}_M$-modules. 
This invokes the notion of a differential graded (dg) vertex algebroid.
The construction of $\CDO_{M,\xi}$ is due to Gorbounov, Malikov and 
Schechtman, and so is the notion of a vertex algebroid;
the dg version is first introduced here.
Section 4 contains the main result, namely the construction of a sheaf of 
dg conformal vertex algebras that provides a fine resolution of 
$\CDO_{M,\xi}$.
This `infinite dimensional Dolbeault complex' plays a role for the Witten
genus similar to that of the Dolbeault complex for the Todd genus.
\end{abstract}

\section{Background \& Overview of the Paper}

In this section, all manifolds are compact unless specified otherwise.

In a study of connections between topology and physics, Witten considered 
various types of conformal field theory associated to a manifold $M$.
Among them, one type
\footnote{
$N=1/2$ sigma model with only right-moving fermions.
}
can only be constructed when $M$ is spin and $p_1(M)=0$.
The torus partition function of this conformal field theory is up to 
a constant factor a modular form of weight $\frac{1}{2}\dim M$. 
The $q$-expansion of the modular form is given by
\beq \label{Witten.genus}
W(M)=\int_M
  \hat{A}(TM)\,
  ch\left(\bigotimes_{k=1}^\infty\t{Sym}_{q^k}(TM\otimes\bb{C})\right)
  \cdot\prod_{k=1}^\infty(1-q^k)^{\dim M}.
\eeq
\cite{Witten}
This is the Witten genus of $M$.
Unfortunately, physicists' construction of the conformal field theory
and their reasoning for the modularity of (\ref{Witten.genus}) involve
the ill-defined path integral.
From a mathematical point of view, the expression in (\ref{Witten.genus})
makes sense whenever $M$ is oriented but is not a priori modular.
Zagier gave a mathematical proof that 
\beq \label{Witten.modular}
W(M)\t{ is a modular form over }\bb{Z}\t{ when }M\t{ is string}.
\eeq
\footnote{
Let $\lambda\in H^4(B\t{Spin};\bb{Z})\cong\bb{Z}$ denote the generator 
twice which is $p_1$, and also the corresponding characteristic class 
for spin vector bundles.
A spin vector bundle or manifold is said to be string if its $\lambda$
vanishes.
Moreover, for $n\geq 1$, there is a topological group $\t{String}_n$ that 
admits a homomorphism into $\t{Spin}_n$ such that the induced map 
$B\t{String}_n\rightarrow B\t{Spin}_n$ is the homotopy fiber of 
$\lambda|_{B\t{Spin}_n}:B\t{Spin}_n\rightarrow K(\bb{Z},4)$.
For specific models of $\t{String}_n$, see \cite{ST,BCSS,Henriques}.
A string structure on a spin vector bundle or manifold is a lifting of 
its structure group from $\t{Spin}_n$ to $\t{String}_n$.
}
\cite{Zagier}
However, a deeper understanding of this result is desired.
Guidance has been provided by an analogy with the $\A$-genus.
Defined for any oriented manifold $M$, $\A(M)$ is a priori a rational 
number but
\beq \label{Ahead.integer}
\A(M)\t{ is an integer when }M\t{ is spin}.
\eeq
It is well known that the search for an explanation of (\ref{Ahead.integer}) 
largely motivated the invention of topological $K$-theory and 
the discovery of the Atiyah-Singer index theorem.

A topological explanation of (\ref{Ahead.integer}) consists of two 
ingredients: 
$KO$-theory and a $KO$-orientation of spin vector bundles.
The orientation can be regarded as a map of ring spectra from $M\t{Spin}$ 
to $KO$.
The induced map in homotopy defines a $\pi_*(KO)$-valued bordism invariant 
for spin manifolds $M$.
If the $2$-torsions in $\pi_*(KO)$ are ignored (or equivalently, if $\dim M$ 
is divisible by $4$), this invariant equals $\A(M)$, proving that it is 
an integer.
A topological explanation of (\ref{Witten.modular}) consists of similar 
ingredients:
the ring spectrum tmf and the $\sigma$-orientation
\beq \label{sigma}
\sigma:M\t{String}\rightarrow\t{tmf}.
\eeq
The construction of tmf and the $\sigma$-orientation is the culmination of 
the combined work of many people.~\cite{LRS,AHS,Hopkins,AHR}
Modularity is built into tmf in the sense that $\pi_*(\t{tmf})$ maps into 
the ring of modular forms and it is an isomorphism away from $6$.
The map induced by (\ref{sigma}) in homotopy defines 
a $\pi_*(\t{tmf})$-valued bordism invariant for string manifolds $M$.
Once $6$ is inverted, this invariant equals $W(M)$, which is therefore
a modular form.~\cite{Hopkins,AHR}

Of course, (\ref{Ahead.integer}) also has a geometric explanation.
It is closely related to the interpretation of $KO$-theory in terms 
of Clifford module bundles.
Namely, the $KO$-orientation of a spin manifold is represented by 
the symbol of its Clifford-linear Dirac operator, and the associated 
$\pi_*(KO)$-valued invariant by the kernel of 
this operator.~\cite{spin.geom}
In contrast, tmf lacks a description similar to that of $K$- or $KO$-theory
using vector bundles.
In light of the physical origin of the Witten genus, a geometric explanation
of (\ref{Witten.modular}) and a related description of tmf seem to entail 
a mathematical definition of conformal field theory.

Vertex algebras provide a mathematical approach to conformal field theory 
\cite{Kac.VA,Frenkel.Ben-Zvi} and may therefore provide a geometric 
meaning of the Witten genus and tmf.
The first result of this kind was due to Gorbounov, Malikov and Schechtman.
Here is an outline of the part of their work in \cite{GMS1} related to 
the Witten genus.
The space of chiral differential operators $\CDO(U)$ over an open set $U$ 
in $\bb{C}^d$ is a basic example of a conformal vertex algebra
(\S\ref{sec.alg.CDO}-\S\ref{sec.holo.CDO}).
To patch these local objects into a sheaf over a complex manifold $M^d$, 
it turns out to require the choice of some geometric data $\xi$
(\S\ref{sec.CDO.mfld}) whose existence requires $c_1(M)=c_2(M)=0$
rationally (\S\ref{sec.obstructions}-\S\ref{sec.char.classes}).
Denote the resulting sheaf of conformal vertex algebras by $\CDO_{M,\xi}$.
Its cohomology $H^*(\CDO_{M,\xi})$ is a conformal vertex superalgebra with 
character (\S\ref{sec.VA})
\beq \label{Witten.genus.geometric}
\t{char}\,H^*(\CDO_{M,\xi})=W(M)\cdot(\t{a constant factor})
\eeq
where the constant factor depends only on $d$ (\S\ref{sec.CDO.mfld.Euler}).
This provides a geometric meaning of the Witten genus for a particular class
of complex manifolds.

The author's present goal is to extend the result 
(\ref{Witten.genus.geometric}) and obtain an interpretation of the Witten 
genus for a larger class of manifolds that includes string manifolds.
To motivate the first step that has been taken, again consider an analogy.
For complex manifolds, the appropriate analogue of the Witten genus is 
the Todd genus.
By the Hirzebruch-Riemann-Roch theorem, the Todd genus of a complex manifold
$M$ has the following closely related interpretations
\begin{align*}
\t{Td}(M) 
&=\t{Euler characteristic of the sheaf of holomorphic functions }\mc{O}_M\\
&=\t{Euler characteristic of the Dolbeault complex }(\E^{0,*}(M),\dbar)\\
&=\t{super dimension of the kernel of the operator }\dbar+\dbar^*\t{ on }
  \E^{0,*}(M);
\end{align*}
\footnote{
The formal adjoint $\dbar^*$ of $\dbar$ is defined with respect to 
a hermitian metric on $M$.
}
among them, it is the last one that has a generalization to all spin$^c$ 
manifolds by the Atiyah-Singer index theorem.~\cite{spin.geom}
In the case of the Witten genus, each sheaf of chiral differential operators
$\CDO_{M,\xi}$ plays the role of $\mc{O}_M$.
In fact, its lowest-weight component is $(\CDO_{M,\xi})_0=\mc{O}_M$
(\S\ref{sec.CDO.iso}, \S\ref{sec.CDO.mfld}).
The analogy with the Todd genus suggests that we find an analogue of 
the Dolbeault complex.
More precisely, for each $\CDO_{M,\xi}$, we would like to construct
\beq \label{desired.CDC}
\begin{array}{ll}
\multicolumn{2}{l}{\t{a sheaf of dg conformal vertex algebras such that}}\\
\t{i.} & \hspace{-0.07in}
\t{it provides a fine resolution of the sheaf of conformal vertex algebras }
\CDO_{M,\xi}\,\t{, and} \\
\t{ii.} & \hspace{-0.07in}
\t{its weight-zero component is the Dolbeault resolution of }\mc{O}_M
\end{array}
\eeq
The construction of such an object, which is outlined below, is 
the main result of this article.

At first glance, there are two difficulties in building the desired object
(\ref{desired.CDC}).
First, we have to deal with an infinite number of sheaves --- 
i.e.~the weight components $(\CDO_{M,\xi})_k$, $k\geq 0$ ---
as well as an infinite number of morphisms between these sheaves given by 
the vertex algebra structure.
Second, all positive weight components are not $\mc{O}_M$-modules
(\S\ref{sec.wt.subsheaves}).
Both difficulties are tackled using vertex algebroids.
The notion of a vertex algebroid was also introduced by Gorbounov, 
Malikov and Schechtman.~\cite{GMS2}
\footnote{
The author apologizes for using a different set of notations.
The structure maps denoted by $\ast,\{\;\}_0,\{\;\}_1$ in this article 
(\S\ref{sec.Void}) are respectively equal to the maps 
$-\gamma$, $\langle\;\rangle$, $-c+\frac{1}{2}\d\circ\langle\;\rangle$
in \cite{GMS2}.
}
Roughly speaking, given a vertex algebra, the part of its structure 
involving only the two lowest weights can be rephrased in terms of 
a ring and two modules (\S\ref{sec.Void}-\S\ref{sec.Void.morphism}).
The resulting algebraic notion is called a vertex algebroid.
The forgetful functor
\beqq
\xymatrix{
  \t{category of vertex algebras}\quad
  \ar[rr]^{\t{forget}}
  &&
  \quad\t{category of vertex algebroids}
  \ar@/^/[ll]^{\t{free}}
}
\eeqq
has a left adjoint that defines the vertex algebra `freely generated'
by a vertex algebroid (\S\ref{sec.VA.VAoid.functors}).
For example, the space of chiral differential operators $\CDO(U)$ over 
$U\subset\bb{C}^d$ is freely generated by a vertex algebroid defined in
terms of $\mc{O}(U)$, $\Omega^1(U)$, $\T(U)$ --- i.e.~the holomorphic 
functions, $1$-forms and vector fields on $U$ (\S\ref{sec.CDO.VAoid}). 
The above discussion has a similar dg version. 
For concreteness, a dg vertex algebroid consists of data of the form
\beqq
\big(A^*,\,\Omega^*,\,\T^*,\,
  \Delta:\T^*\rightarrow\Omega^{*+1},\,
  \ast:A^*\otimes\T^*\rightarrow\Omega^*,\,
  \{\;\}_0:\T^*\otimes\T^*\rightarrow A^*,\,
  \{\;\}_1:\T^*\otimes\T^*\rightarrow\Omega^*\big)
\eeqq
satisfying a host of axioms.
In particular, $A^*$ is a dg ring, $\Omega^*,\T^*$ are $A^*$-modules,
$\Delta$ is a chain map of degree one, and $\ast,\{\;\}_0,\{\;\}_1$ 
are null homotopies of chain maps constructed from $\Delta$ 
(\S\ref{sec.dgVAoid1}).

The construction of the desired object (\ref{desired.CDC}) is summarized 
in figure \ref{fig.CDC.construction}.
\footnote{
Notations: 
$\check C^*(-)$ are \v Cech complexes of a finite good open cover;
$\mc{O}_M$, $\Omega^1_M$, $\T_M$ are the sheaves of holomorphic functions,
$1$-forms and vector fields;
$\E^{0,*}_M$ is the Dolbeault resolution of $\mc{O}_M$;
tensor products are taken over $\mc{O}_M$.
}
The diagram starts with the dg vertex algebra $\check C^*(\CDO_{M,\xi})$.
\footnote{
In fact, the \v Cech complex valued in a sheaf of vertex algebras is only 
a dg vertex algebra up to homotopy (\S\ref{sec.dgVA}).
Similarly, every term in the right column of figure 
\ref{fig.CDC.construction} except the last is only a dg vertex algebroid 
up to homotopy.
}
The top row follows from the comment about $\CDO(U)$ in the previous 
paragraph (\S\ref{sec.CDO.dgVAoid.Cech}).
The right column is a sequence of quasi-isomorphisms of dg vertex
algebroids (\S\ref{sec.dgVAoid.hom}).
The first quasi-isomorphism extends 
$\check\Delta,\check\ast,\{\;\check\}_0,\{\;\check\}_1$ from \v Cech 
complexes to \v Cech-Dolbeault complexes such that $\check\Delta$ has 
\v Cech degree one and Dolbeault degree zero 
(\S\ref{sec.dgVAoid.Cech.Dolbeault.1}).
The second quasi-isomorphism modifies $\check\Delta$ by a homotopy into 
a map $\bar\Delta$ with \v Cech degree zero and Dolbeault degree one; 
the other three maps are also modified appropriately 
(\S\ref{sec.dgVAoid.Cech.Dolbeault.2}).
This step requires the choice of some geometric data $(\nabla,H)$ 
to be explained below.
The third quasi-isomorphism restricts 
$\bar\Delta,\bar\ast,\{\;\bar\}_0,\{\;\bar\}_1$ to global sections 
(\S\ref{sec.dgVAoid.Dolbeault}).
Finally, the bottom arrow defines a new dg vertex algebra 
$\Gamma(\CDC_{M,\nabla,H})$.
As the notation suggests, it is in fact the global sections of 
a sheaf of dg vertex algebras $\CDC_{M,\nabla,H}$ 
(\S\ref{sec.dgVAoid.Dolbeault}-\S\ref{sec.dgVA.Dolbeault}).
This object has the properties stated in (\ref{desired.CDC}).
\footnote{
The steps taken to construct this object may seem devious.
However, even though the presheaf of vertex algebroids induced by 
$\CDO_{M,\xi}$ is defined in terms of $\mc{O}_M$-modules, it contains 
collections of maps that do not form morphisms of sheaves 
(\S\ref{sec.presheaf.VAoid}) and hence do not extend easily to maps 
between Dolbeault resolutions.
Instead, those data define non-chain maps between \v Cech complexes, namely 
$\check\ast,\{\;\check\}_0,\{\;\check\}_1$ in the first dg vertex 
algebroid in figure \ref{fig.CDC.construction} 
(\S\ref{sec.CDO.dgVAoid.Cech}).
Their extensions to \v Cech-Dolbeault complexes respect the Dolbeault 
but not the \v Cech differential (\S\ref{sec.dgVAoid.Cech.Dolbeault.1}).
They are then modified into maps $\bar\ast,\{\;\bar\}_0,\{\;\bar\}_1$ 
that respect the \v Cech but not the Dolbeault differential 
(\S\ref{sec.dgVAoid.Cech.Dolbeault.2}).
In particular, $\bar\ast,\{\;\bar\}_0,\{\;\bar\}_1$ define morphisms of 
sheaves between Dolbeault resolutions (\S\ref{sec.dgVAoid.Dolbeault}).
}
In particular, there is a fine resolution
\beq \label{CDO.resolution}
\xymatrix{
  0\ar[r] &
  \CDO_{M,\xi}\ar[r]^{i\phantom{a}} &
  \CDC_{M,\nabla,H}
}
\eeq
where $i$ is a morphism of sheaves of conformal vertex superalgebras 
(\S\ref{sec.dgVAoid.Dolbeault}-\S\ref{sec.dgVA.Dolbeault}, 
\S\ref{sec.CDC.conformal}).

\begin{figure}[t] 
{\fn
\beqq
\xymatrix{
\t{\underline{dg vertex algebras}} && 
\t{\underline{dg vertex algebroids}} \\ 
\hspace{0.3in}\check C^*(\CDO_{M,\xi})\hspace{0.3in}
  \ar[rr]^{\t{forget}\phantom{aaaaaaaaaaaaaa}} &&
\quad\big(
  \check C^*(\mc{O}_M),\,
  \check C^*(\Omega^1_M),\,
  \check C^*(\T_M),\,
  \check\Delta,\,\check\ast,\,\{\;\check\}_0,\,\{\;\check\}_1
\big)\quad
  \ar[d]^\simeq_\cap\ar@/^/[ll]^{\t{free}\phantom{aaaaaaaaaaaaaaa}} \\
&&
\big(
  \check C^*(\E^{0,*}_M),\,
  \check C^*(\Omega^1_M\otimes\E^{0,*}_M),\,
  \check C^*(\T_M\otimes\E^{0,*}_M),\,
  \check\Delta,\,\check\ast,\,\{\;\check\}_0,\,\{\;\check\}_1
\big)\ar[d]^\simeq \\
&&
\big(
  \check C^*(\E^{0,*}_M),\,
  \check C^*(\Omega^1_M\otimes\E^{0,*}_M),\,
  \check C^*(\T_M\otimes\E^{0,*}_M),\,
  \bar\Delta,\,\bar\ast,\,\{\;\bar\}_0,\,\{\;\bar\}_1
\big) \\
\hspace{0.15in}\Gamma(\CDC_{M,\nabla,H})\hspace{0.15in}
&&
\quad\big(
  \Gamma(\E^{0,*}_M),\,
  \Gamma(\Omega^1_M\otimes\E^{0,*}_M),\,
  \Gamma(\T_M\otimes\E^{0,*}_M),\,
  \bar\Delta,\,\bar\ast,\,\{\;\bar\}_0,\,\{\;\bar\}_1
\big)\quad
  \ar[u]^\cup_\simeq\ar[ll]^{\t{free}\phantom{aaaaaaaaaaaaaaaaaaa}}
}
\eeqq}
\caption{
  Construction of the chiral Dolbeault complex.
  \label{fig.CDC.construction}}
\end{figure}

The following are some features in the construction of 
$\CDC_{M,\nabla,H}$.
{\bf (i)} The input data include a connection $\nabla$ on $TM$ and 
a $3$-form $H$ satisfying certain conditions, most notably 
$dH=-\Tr(R\wedge R)$, where $R$ is the curvature of $\nabla$ 
(\S\ref{sec.h}-\S\ref{sec.B}).
A choice of $(\nabla,H)$ is essentially equivalent to a choice of the data 
$\xi$ needed to construct $\CDO_{M,\xi}$ (\S\ref{sec.CDC.iso.classes}).
For $(\nabla,H)$ to exist, a necessary condition is $c_1(M)=c_2(M)=0$ 
rationally, which is also sufficient in the case $M$ is K\"ahler 
(\S\ref{sec.CDC.Kahler}).
{\bf (ii)} While the definition of each $\CDO_{M,\xi}$ depends 
necessarily on local coordinates (\S\ref{sec.CDO.mfld}), 
$\CDC_{M,\nabla,H}$ has a global description
(\S\ref{sec.dgVAoid.Dolbeault}).
{\bf (iii)}
The differentials on the weight-zero and -one components are respectively
the ordinary Dolbeault operator on $\E^{0,*}_M$ and a `deformed' Dolbeault
operator (\S\ref{sec.dgVA.Dolbeault})
\beqq
\dbardef\;\t{on}\;\;(\Omega^1_M\oplus\T_M)\otimes\E^{0,*}_M,
\qquad
(\alpha,X)\mapsto(\dbar\alpha+\bar\Delta X,\dbar X).
\eeqq
{\bf (iv)} It follows from the resolution (\ref{CDO.resolution}) that 
$\Gamma(\CDC_{M,\nabla,H})$ computes $H^*(\CDO_{M,\xi})$.
This provides a new interpretation of the conformal vertex superalgebra
$H^*(\CDO_{M,\xi})$ and hence, by (\ref{Witten.genus.geometric}), 
the Witten genus
\beq \label{CDC.Witten}
\t{char}\,H^*\big(\Gamma(\CDC_{M,\nabla,H})\big)
=W(M)\cdot(\t{a constant factor}).
\eeq
In view of the connection between conformal vertex (super)algebras and 
conformal field theory, (\ref{CDC.Witten}) goes some way to `explaining'
the modularity of $W(M)$.

This paragraph is added in February 2010.
After completing the work reported here, the author has realized that 
the construction of $\CDC_{M,\nabla,H}$ is somewhat ad hoc and it is not 
clear how to generalize the interpretation of the Witten genus in
(\ref{CDC.Witten}).
Another construction $\CDCnew_{M,\nabla,H}$ that also fulfils 
(\ref{desired.CDC}) is given in \cite{newpaper}.
Roughly speaking, it is a sheaf of `smooth' chiral differential operators
on the cs-manifold $\Pi\,\overline{TM}$, equipped with a derivation
associated to an odd vector field.
While $\CDC_{M,\nabla,H}$ embeds quasi-isomorphically into 
$\CDCnew_{M,\nabla,H}$, the latter belongs to a more general and 
systematic framework.
In fact, as a sheaf of conformal vertex superalgebras 
$\CDCnew_{M,\nabla,H}$ is defined whenever $M$ is almost complex and
$c_1(M)=c_2(M)=0$ rationally (but the derivation may not be order-two
in general).
On the other hand, some of the ideas in \cite{newpaper} came from
this work.
For a discussion of future directions, see also \cite{newpaper}.

The author would like to thank Matthew Ando, Victor Ka\v c, Haynes Miller,
Stephan Stolz, Peter Teichner and his Ph.D. advisor Ralph Cohen for their 
useful comments as well as their continual interest and encouragement 
in this project.
The notion of a vertex algebroid was first brought to the author's 
attention by Haynes Miller.
The author worked out most of the content of this article when he was 
a Moore Instructor at MIT before completing it at the Max-Planck-Institut 
f\"ur Mathematik;
he is grateful for the generous support from both institutes.

\setcounter{equation}{0}
\section{Chiral Differential Operators}

This section reviews the construction of chiral differential operators 
on complex manifolds as well as their relation with the Witten genus.
The presentation here somewhat differs from \cite{GMS1}.

\subsection{} \label{sec.VA}
In this article, we adopt the definition of a \emph{vertex superalgebra} 
and a \emph{conformal vertex superalgebra} in \cite{Frenkel.Ben-Zvi}.
This in particular means the state space of a vertex algebra always has 
a $\bb{Z}_+$-grading, referred to as \emph{weight}.
(See also \cite{Kac.VA}.)
The \emph{character} of a conformal vertex superalgebra $V$, when defined, 
is the formal power series
\beqq
\t{char}\;V=\sum_k q^{k-c/24}(\dim V_k^{\t{even}}-\dim V_k^{\t{odd}})
  \in q^{-c/24}\cdot\bb{Z}[[q]]
\eeqq
where $V^{\t{even}}_k,V^{\t{odd}}_k$ are the even and odd parts of 
the weight-$k$ component, and $c$ is the central charge.

\subsection{} \label{sec.alg.CDO}
Let $\mc{W}$ be the unital associative $\bb{C}$-algebra with generators 
$a_{in}, b^i_n$ for $n\in\bb{Z}$, $i=1,\ldots,d$, and relations
\beqq
[a_{in},b^j_m]=\delta^j_i\delta_{n,-m},\qquad
[a_{in},a_{jm}]=0=[b^i_n,b^j_m].
\eeqq
The subalgebra $\mc{W}_+\subset\mc{W}$ generated by $a_{in}$, $n\geq 0$, 
and $b^i_n$, $n>0$, is commutative and hence admits a trivial 
representation $\bb{C}$. 
The induced representation $\mc{W}\otimes_{\mc{W}_+}\bb{C}$ of $\mc{W}$ 
has the structure of a conformal vertex algebra. 
Let us denote it by $\CDO_{\t{alg}}(\bb{C}^d)$ and 
define its structure below.

A general element of $\CDO_{\t{alg}}(\bb{C}^d)$ is a linear combination of 
elements of the form
\beq \label{general.element}
a_{i_1 n_1}\cdots a_{i_k n_k}b^{j_1}_{m_1}\cdots b^{j_\ell}_{m_\ell}f
\eeq
where $n_1\leq\cdots\leq n_k<0$, $m_1\leq\cdots\leq m_\ell<0$ 
($k,\ell\geq 0$) and $f\in\bb{C}[b^1_0,\ldots,b^d_0]$.
The vacuum is $1$. 
The vertex operators of the elements $a_{i,-1}$ and $b^i_0$ are 
respectively
\beqq
a_i(z)=\sum_{n\in\bb{Z}}a_{in}z^{-n-1}\quad\t{and}\quad
b^i(z)=\sum_{n\in\bb{Z}}b^i_n z^{-n}.
\eeqq
The axioms of a vertex algebra then determine all other vertex operators 
as well as the infinitesimal translation operator. 
The conformal element is
\beq \label{conformal}
a_{1,-1}b^1_{-1}+\ldots+a_{d,-1}b^d_{-1}
\eeq
with central charge $2d$.
The weight zero subspace is 
\beqq
\CDO_{\t{alg}}(\bb{C}^d)_0
=\bb{C}[b^1_0,\ldots,b^d_0]
=\{\,\t{algebraic functions on }\bb{C}^d\} 
\eeqq
and the operators $a_{in}$ and $b^i_n$ change weight by $-n$. 
In particular, the element (\ref{general.element}) has weight 
$-n_1-\ldots-n_k-m_1-\ldots-m_\ell$.

\subsection{} \label{sec.holo.CDO}
According to \S\ref{sec.CDO.VAoid}, there are `holomorphic analogues'
of $\CDO_{\t{alg}}(\bb{C}^d)$.
More precisely, for each open set $U\subset\bb{C}^d$, there is 
a conformal vertex algebra $\CDO(U)$ whose weight zero subspace is
\beqq
\CDO(U)_0=\mc{O}(U)=\{\,\t{holomorphic functions on }U\}
\eeqq
and whose elements are linear combinations of those of the form 
(\ref{general.element}) but with $f\in\mc{O}(U)$. 
This defines a sheaf of conformal vertex algebras $\CDO$ over 
$\bb{C}^d$.

\subsection{} \label{sec.CDO.iso}
Let $\varphi=(\varphi^1,\ldots,\varphi^d):U\rightarrow V$ be 
a biholomorphism between two open sets in $\bb{C}^d$ and 
$g_\varphi$ its matrix derivative, 
i.e.~$(g_\varphi)^i_{\phantom{i}j}=\d_j\varphi^i$.
Define a matrix-valued $1$-form and a (scalar) $3$-form as follows
\beqq
\theta_\varphi=g_\varphi^{-1}\cdot\d g_\varphi,\qquad
WZ_\varphi=\frac{1}{3}\Tr(\theta_\varphi\wedge\theta_\varphi\wedge
  \theta_\varphi).
\eeqq
Notice that $\d\theta_\varphi=-\theta_\varphi\wedge\theta_\varphi$ and
$\d WZ_\varphi=0$.

Consider the problem of finding all isomorphisms of conformal vertex 
algebras $\CDO(V)\rightarrow\CDO(U)$ whose weight-zero component is 
$\varphi^*:\mc{O}(V)\rightarrow\mc{O}(U)$, i.e. $b^i_0\mapsto\varphi^i$.
Such an isomorphism is determined by the images of the weight-one 
elements $a_{i,-1}$. 
The result is stated below.
Firstly, we obtain an isomorphism of vertex algebras (not necessarily 
preserving the conformal elements) if and only if 
$a_{i,-1}$ are sent to
\beq \label{a-1.transform}
a_{j,-1}\big(g_\varphi^{-1}\big)^j_{\phantom{i}i}
+\frac{1}{2}b^j_{-1}
  \Big\{\xi_{jk}+\Tr\big[\theta_\varphi(\d_j)\cdot\theta_\varphi(\d_k)\big]
    \Big\}
  \big(g_\varphi^{-1}\big)^k_{\phantom{i}i}
\eeq
where $\xi_{ij}\in\mc{O}(U)$, $\xi_{ji}=-\xi_{ij}$ and the holomorphic 
$2$-form $\xi=\sum_{i<j}\xi_{ij}db^i_0\wedge db^j_0$ satisfies the equation
\beqq
\d\xi=WZ_\varphi.
\eeqq
According to the Poincar\'e lemma for the $\d$-operator, since 
$\d WZ_\varphi=0$, such $\xi$ exists after possibly replacing $U$ and $V$
by smaller open sets. 
Label this isomorphism by 
\beqq
\varphi^*_\xi:\CDO(V)\rightarrow\CDO(U).
\eeqq
Secondly, $\varphi^*_\xi$ preserves the conformal element 
(\ref{conformal}) if and only if
\beqq
\Tr\theta_\varphi=0.
\eeqq

\subsection{} 
Given isomorphisms $(\varphi_1)^*_{\xi_1}:\CDO(V)\rightarrow\CDO(U)$
and $(\varphi_2)^*_{\xi_2}:\CDO(W)\rightarrow\CDO(V)$ of vertex algebras,
they compose as follows
\beq \label{CDO.iso.comp}
(\varphi_1)^*_{\xi_1}\circ(\varphi_2)^*_{\xi_2}
=(\varphi_2\varphi_1)^*_\eta,\qquad
\eta=\xi_1+\varphi_1^*\xi_2+\sigma_{\varphi_2,\varphi_1}
\eeq
where $\sigma_{\varphi_2,\varphi_1}=
\Tr\big(\theta_{\varphi_1}\wedge g_{\varphi_1}^{-1}\cdot
\varphi_1^*\theta_{\varphi_2}\cdot g_{\varphi_1}\big)$ is a $2$-form on $U$.

\subsection{} \label{sec.CDO.mfld}
Let $M^d$ be a compact, complex manifold. 
Choose a finite good (in the sense of \cite{Bott.Tu}) open cover 
$\mf{U}=\{U_1,\ldots,U_N\}$ of $M$ with holomorphic coordinate charts 
$\varphi_\alpha:U_\alpha\rightarrow\bb{C}^d$, $\alpha=1,\ldots,N$. 
Let 
\beqq
W_\alpha=\varphi_\alpha(U_\alpha),\quad
W_{\alpha\beta}=\varphi_\alpha(U_\alpha\cap U_\beta)\quad\t{and}\quad 
W_{\alpha\beta\gamma}=\varphi_\alpha(U_\alpha\cap U_\beta\cap U_\gamma)
\eeqq 
whenever $U_\alpha\cap U_\beta\neq\varnothing$ and 
$U_\alpha\cap U_\beta\cap U_\gamma\neq\varnothing$.
For the coordinate transformations 
\beqq
\varphi_{\beta\alpha}
=\varphi_\beta\circ\varphi_\alpha^{-1}|_{W_{\alpha\beta}}:
  W_{\alpha\beta}\rightarrow W_{\beta\alpha}
\eeqq
we shall write $g_{\beta\alpha}$, $\theta_{\beta\alpha}$, 
$WZ_{\beta\alpha}$, $\sigma_{\gamma\beta\alpha}$ respectively for 
$g_{\varphi_{\beta\alpha}}$, $\theta_{\varphi_{\beta\alpha}}$, 
$WZ_{\varphi_{\beta\alpha}}$, 
$\sigma_{\varphi_{\gamma\beta},\varphi_{\beta\alpha}}$.

Assign to each $U_\alpha$ the sheaf $\CDO|_{W_\alpha}$.
To glue these into a sheaf of conformal vertex algebras over $M$, it 
follows from \S\ref{sec.CDO.iso} that we need to choose 
on each $W_{\alpha\beta}$ a holomorphic solution $\xi_{\beta\alpha}$ to
\beq \label{d.2form.WZ}
d\xi_{\beta\alpha}=\d\xi_{\beta\alpha}=WZ_{\beta\alpha}
\eeq
in order to define an isomorphism 
$(\varphi_{\beta\alpha})^*_{\xi_{\beta\alpha}}:\CDO(W_{\beta\alpha})
\rightarrow\CDO(W_{\alpha\beta})$ of vertex algebras.
This is possible, perhaps after a refinement of $\mf{U}$.
The following conditions must be satisfied:
\beq
\label{glue.assoc}
& \varphi_{\beta\alpha}^*\xi_{\gamma\beta}-\xi_{\gamma\alpha}
  +\xi_{\beta\alpha}+\sigma_{\gamma\beta\alpha}=0 &
  \quad\t{on each }W_{\alpha\beta\gamma} \\
\label{glue.conformal}
& \Tr\theta_{\beta\alpha}=0 &
  \quad\t{on each }W_{\alpha\beta}
\eeq
Indeed, (\ref{glue.assoc}) is equivalent to the cocycle condition 
$(\varphi_{\gamma\alpha})^*_{\xi_{\gamma\alpha}}=
(\varphi_{\beta\alpha})^*_{\xi_{\beta\alpha}}\circ
(\varphi_{\gamma\beta})^*_{\xi_{\gamma\beta}}$ in view of 
(\ref{CDO.iso.comp}), and (\ref{glue.conformal}) ensures that the conformal
elements from $U_\alpha$ agree over intersections, 
according to the last statement in \S\ref{sec.CDO.iso}.
These conditions allow us to construct a sheaf of conformal vertex 
algebras on $M$ that restricts to $\CDO|_{W_\alpha}$ on each $U_\alpha$. 
Denote this sheaf by $\CDO_{M,\xi}$ where $\xi$ stands for the \v Cech 
cochain $\{\xi_{\beta\alpha}\}$. 
\footnote{
For each open set $U\subset M$, an element of $\CDO_{M,\xi}(U)$ is 
a collection of elements of $\CDO(\varphi_\alpha(U\cap U_\alpha))$ that 
agree over intersections via 
$(\varphi_{\beta\alpha})^*_{\xi_{\beta\alpha}}$.
The finiteness of $\mf{U}$ ensures that $\CDO_{M,\xi}(U)$ satisfies 
all the axioms of a vertex algebra.
}

\subsection{} \label{sec.obstructions}
Denote by $\Omega_M^p$ the sheaf of holomorphic $p$-forms on $M$ and 
$\Omega_M^{p,\t{cl}}$ its subsheaf of closed forms.
Suppose $\xi=\{\xi_{\beta\alpha}\}$ satisfying (\ref{d.2form.WZ}) have 
been chosen but (\ref{glue.assoc})-(\ref{glue.conformal}) do not 
necessarily hold.

Consider the \v Cech cochain
\beq \label{cocycle.assoc}
\big\{
  \varphi_{\beta\alpha}^*\xi_{\gamma\beta}-\xi_{\gamma\alpha}
  +\xi_{\beta\alpha}+\sigma_{\gamma\beta\alpha}
\big\}_{\alpha\beta\gamma}
  \in\check C^2(\mf{U},\Omega_M^{2,\t{cl}}).
\eeq
Indeed, (\ref{cocycle.assoc}) consists of closed forms in view of 
(\ref{d.2form.WZ}) and (\ref{CDO.iso.comp}).
It follows from the associativity of (\ref{CDO.iso.comp}) that 
(\ref{cocycle.assoc}) is a cocycle.
Keeping the coordinate charts $\varphi_\alpha$ fixed, one can obtain 
the equations in (\ref{glue.assoc}) after a modification of $\xi$ 
if and only if (\ref{cocycle.assoc}) is a coboundary 
of closed $2$-forms.
Therefore the obstruction to constructing a sheaf of (not necessarily 
conformal) vertex algebras $\CDO_{M,\xi}$ as described in 
\S\ref{sec.CDO.mfld} is the class of (\ref{cocycle.assoc}) in 
$H^2(\Omega_M^{2,\t{cl}})$.
This class does not depend on $\xi$.

Now consider the \v Cech cochain
\beq \label{cocycle.conformal}
\big\{\Tr\theta_{\beta\alpha}\big\}_{\alpha\beta}
  \in\check C^1(\mf{U},\Omega^{1,\t{cl}}_M).
\eeq
It is easy to verify that (\ref{cocycle.conformal}) indeed consists of 
closed forms and is a cocycle.
Furthermore, one can obtain the equations in (\ref{glue.conformal}) after 
a modification of the coordinate charts $\varphi_\alpha$ if and only if 
(\ref{cocycle.conformal}) is a coboundary of closed $1$-forms. 
Therefore given a sheaf of vertex algebras $\CDO_{M,\xi}$ as constructed 
above, its obstruction to having a global conformal element is the 
class of (\ref{cocycle.conformal}) in $H^1(\Omega_M^{1,\t{cl}})$.

\subsection{} \label{sec.char.classes} 
Denote by $\E_M^p$ the sheaf of smooth $p$-forms on $M$ and 
$\E_M^{p,\t{cl}}$ its subsheaf of closed forms.
Choose a connection on $TM$ and let $\Gamma_\alpha$ be the connection 
$1$-form on $W_\alpha$ associated to the coordinate chart $\varphi_a$.
The curvature operator is locally given by
$R_\alpha=d\Gamma_\alpha+\Gamma_\alpha\wedge\Gamma_\alpha$.
Define a $3$-form on $W_\alpha$ by
\beqq
CS(\Gamma_\alpha)=\Tr(\Gamma_\alpha\wedge R_\alpha)
  -\frac{1}{3}\Tr(\Gamma_\alpha\wedge\Gamma_\alpha\wedge\Gamma_\alpha)
\eeqq

Consider the computation in the \v Cech-de Rham double complex 
$\check C^*(\mf{U},\E_M^*)$ shown in figure \ref{fig.CdR.ch2}:
\footnote{
The total differential of the \v Cech-de Rham complex is 
$\delta+(-1)^p d$, where $\delta$ is the \v Cech differential, 
$d$ the de Rham differential and $p$ the \v Cech degree.
}
it shows that (\ref{cocycle.assoc}) as an element of 
$\check C^2(\mf{U},\E_M^2)$ is cohomologous to 
$\{-\Tr(R_\alpha\wedge R_\alpha)\}_\alpha\in\check C^0(\mf{U},\E_M^4)$
in the total \v Cech-de Rham complex.
By Chern-Weil theory, the latter represents $8\pi^2 ch_2(M)$.
In other words, the obstruction class represented by (\ref{cocycle.assoc}) 
maps to $8\pi^2 ch_2(M)$ under the following composition
\beqq
\xymatrix{
  H^2(\Omega_M^{2,\t{cl}})\ar[r] &
  H^2(\E_M^{2,\t{cl}})\ar[r]^\cong &
  H^4(M;\bb{C})
}
\eeqq
where the first arrow is induced by an inclusion of sheaves and 
the second by the inclusion of the \v Cech complex valued in 
$\E_M^{2,\t{cl}}$ into the total \v Cech-de Rham complex.

\begin{figure}[t] 
{\scriptsize
\beqq
\begin{array}{r|ccccc}
4 & 
\big\{\,\Tr(R_\alpha\wedge R_\alpha)\,\big\}_\alpha &&&& 
(\uparrow \t{de Rham};\; \rightarrow \t{\v Cech}) \\
& \uparrow_d &&&& \\
3 &
\big\{\,CS(\Gamma_\alpha)\,\big\}_\alpha & 
\stackrel{\delta}{\longrightarrow} &
\big\{\,WZ_{\beta\alpha}+d\,\Tr(\theta_{\beta\alpha}\wedge\Gamma_\alpha)
  \,\big\}_{\alpha\beta} && \\
&&& \uparrow_d && \\
2 &&&
\big\{\,\xi_{\beta\alpha}
  +\Tr(\theta_{\beta\alpha}\wedge\Gamma_\alpha)\,\big\}_{\alpha\beta} &
\stackrel{\delta}{\longrightarrow} &
\big\{
  \varphi^*_{\beta\alpha}\xi_{\gamma\beta}
  -\xi_{\gamma\alpha}
  +\xi_{\beta\alpha}
  -\sigma_{\gamma\beta\alpha}
\big\}_{\alpha\beta\gamma} \\ 
\vdots &&&&& \\
\hline 
& 0 && 1 && 2
\end{array}
\eeqq}
\caption{
  \v Cech-de Rham cochains related to the construction of $\CDO_{M,\xi}$.
  \label{fig.CdR.ch2}}
\end{figure}

Similarly, the computation in figure \ref{fig.CdR.ch1} shows that 
(\ref{cocycle.conformal}) as an element of $\check C^1(\mf{U},\E^1_M)$ 
is cohomologous to $\{-\Tr R_\alpha\}_\alpha\in\check C^0(\mf{U},\E_M^2)$.
The latter represents $2\pi i ch_1(M)$.
In other words, the obstruction class represented by 
(\ref{cocycle.conformal}) maps to $2\pi i ch_1(M)$ under the following 
composition
\beqq
\xymatrix{
  H^1(\Omega_M^{1,\t{cl}})\ar[r] &
  H^1(\E_M^{1,\t{cl}})\ar[r]^\cong &
  H^2(M;\bb{C})
}
\eeqq
where the first arrow is again induced by an inclusion of sheaves and 
the second by the inclusion of the \v Cech complex valued in
$\E_M^{1,\t{cl}}$ into the total \v Cech-de Rham complex.

\begin{figure}[t] 
{\scriptsize
\beqq
\begin{array}{r|ccc}
2 &
\{\Tr R_\alpha\}_\alpha && \\
& \uparrow_d && \\
1 &
\{\Tr\Gamma_\alpha\}_\alpha &
\stackrel{\delta}{\longrightarrow} &
\big\{\Tr\theta_{\beta\alpha}\big\}_{\alpha\beta} \\
\vdots &&& \\
\hline
& 0 && 1 
\end{array}
\eeqq}
\caption{
  \v Cech-de Rham cochains related to the conformal structure of 
  $\CDO_{M,\xi}$.
  \label{fig.CdR.ch1}}
\end{figure}

\subsection{} \label{sec.CDO.iso.classes}
Assume that both obstructions discussed in \S\ref{sec.obstructions} 
vanish.
The cochains $\xi$ satisfying (\ref{d.2form.WZ}) and (\ref{glue.assoc}) 
form an affine space modeled on the space of cocycles in 
$\check C^1(\mf{U},\Omega_M^{2,\t{cl}})$.
Given two such cochains $\xi,\xi'$, the sheaves of conformal vertex algebras 
$\CDO_{M,\xi},\CDO_{M,\xi'}$ are isomorphic over $M$ if and only if 
there exist closed holomorphic $2$-forms $\xi_\alpha$ on $W_\alpha$ such 
that all diagrams of the form
\beqq
\xymatrix{
  \CDO(W_{\beta\alpha})\ar[rr]^{\t{id}^*_{\xi_\beta}}
    \ar[d]_{(\varphi_{\beta\alpha})^*_{\xi_{\beta\alpha}}} &&
  \CDO(W_{\beta\alpha})
    \ar[d]^{(\varphi_{\beta\alpha})^*_{\xi'_{\beta\alpha}}} \\
  \CDO(W_{\alpha\beta})\ar[rr]^{\t{id}^*_{\xi_\alpha}} &&
  \CDO(W_{\alpha\beta})
}
\eeqq
commute.
By (\ref{CDO.iso.comp}), this is equivalent to the equations
\beqq
\xi'_{\beta\alpha}-\xi_{\beta\alpha}
  =\xi_\alpha-\varphi_{\beta\alpha}^*\xi_\beta.
\eeqq
Namely, $\xi,\xi'$ differ by a coboundary in 
$\check C^1(\mf{U},\Omega_M^{2,\t{cl}})$.
Therefore the isomorphism classes of sheaves of conformal vertex algebras 
$\CDO_{M,\xi}$ constructed as in \S\ref{sec.CDO.mfld} form an affine space 
modeled on $H^1(\Omega_M^{2,\t{cl}})$.
This does not depend on the open cover $\mf{U}$ or the coordinate charts 
$\varphi_\alpha$.

\subsection{} \label{sec.wt.subsheaves}
Let $(\CDO_{M,\xi})_k$ be the weight-$k$ subsheaf of $\CDO_{M,\xi}$, 
$k\geq 0$.
At weight zero, we have
\beq \label{wt0}
(\CDO_{M,\xi})_0=\mc{O}_M.
\eeq
Denote by $\T_M$ the sheaf of holomorphic vector fields on $M$.
There is a short exact sequence
\beq \label{wt1.ses}
\xymatrix{
  0\ar[r] &
  \Omega^1_M\ar[r] &
  (\CDO_{M,\xi})_1\ar[r] &
  \T_M\ar[r] &
  0
}
\eeq
The inclusion of $\Omega^1_M$ sends the $1$-forms $d\varphi_\alpha^i$ 
on $U_\alpha$ to $b_{-1}^i\in\CDO(W_\alpha)=\CDO_{M,\xi}(U_\alpha)$.
The projection to $\T_M$ sends $a_{i,-1}\in\CDO_{M,\xi}(U_\alpha)$ to
the vector fields $\d/\d\varphi_\alpha^i$ on $U_\alpha$.
The way $a_{i,-1}$ transform between different coordinate neighborhoods, 
as described in (\ref{a-1.transform}), indicates that the short exact 
sequence (\ref{wt1.ses}) does not split and $(\CDO_{M,\xi})_1$ is not 
an $\mc{O}_M$-module.
In general, each $(\CDO_{M,\xi})_k$ admits a filtration whose associated 
graded space is an $\mc{O}_M$-module.
\footnote{
The weight-$k$ monomials in $a_{i,-1}$, $b_{-1}^i$ admit a partial 
ordering with the property that each isomorphism $\varphi^*_\xi$ sends 
a monomial to a sum of monomials of the same or lower order.
Therefore the weight-$k$ monomials from all $\CDO_{M,\xi}(U_\alpha)$ 
equal to or lower than a particular order span a subsheaf. 
These subsheaves form the described filtration of $(\CDO_{M,\xi})_k$.
}
For example, $(\CDO_{M,\xi})_2$ can be filtered as follows, with 
the indicated successive quotients
\beq \label{wt2.filtration}
\xymatrix{
  0\ar[rr]^{\t{Sym}^2\Omega^1_M} &&
  \mc{F}^1\ar[rr]^{\Omega^1_M} &&
  \mc{F}^2\ar[rr]^{\T_M\otimes\Omega^1_M} &&
  \mc{F}^3\ar[rr]^{\T_M} &&
  \mc{F}^4\ar[rr]^{\t{Sym}^2\T_M} &&
  (\CDO_{M,\xi})_2
}
\eeq

\subsection{} \label{sec.CDO.mfld.Euler}
The cohomology of a sheaf of conformal vertex algebras such as 
$\CDO_{M,\xi}$ is a conformal vertex superalgebra with the same 
central charge.
The character of $H^*(\CDO_{M,\xi})$ is given by 
\beqq
\t{char}\;H^*(\CDO_{M,\xi})
= q^{-d/12}\sum_{k=0}^\infty q^k\chi\big((\CDO_{M,\xi})_k\big)
\eeqq
where $\chi(-)$ denotes the Euler characteristic.
Because of (\ref{wt0}) and (\ref{wt1.ses}), we have 
\beqq
\chi\big((\CDO_{M,\xi})_0\big)=\chi(\mc{O}_M),\qquad
\chi\big((\CDO_{M,\xi})_1\big)=\chi(\Omega^1_M)+\chi(\T_M).
\eeqq
In general, filtrations like (\ref{wt2.filtration}) allow us to express 
$\chi\big((\CDO_{M,\xi})_k\big)$ in terms of the Euler characteristics 
of $\mc{O}_M$-modules.
Using the Hirzebruch-Riemann-Roch formula, we obtain
\beqq
\t{char}\,H^*(\CDO_{M,\xi})=
q^{-d/12}\chi\left(\bigotimes_{\ell=1}^\infty
  \t{Sym}_{q^\ell}(\Omega^1_M\oplus\T_M)
\right)
=\int_M e^{\frac{1}{2}ch_1(M)}\frac{W(TM)}{\eta(q)^{2d}}
\eeqq
where $\t{Sym}_t(\cdot)=\sum_{n=0}^\infty t^n\t{Sym}^n(\cdot)$ and 
$\eta(q)$ is the Dedekind $\eta$-function.
According to \S\ref{sec.char.classes}, $ch_1(M)=0$ and hence
\beq \label{CDO.mfld.Euler}
\t{char}\;H^*(\CDO_{M,\xi})=\frac{W(M)}{\eta(q)^{2d}}\,.
\eeq
To summarize, given a compact complex manifold $M$, if the obstruction 
classes described in \S\ref{sec.obstructions} vanish, each sheaf of 
conformal vertex algebras $\CDO_{M,\xi}$ provides a geometric 
interpretation of the Witten genus $W(M)$.

\setcounter{equation}{0}
\section{Differential Graded Vertex Algebroids}

The sheaf of vertex algebras $\CDO_{M,\xi}$ is not an $\mc{O}_M$-module 
(\S\ref{sec.wt.subsheaves}), but in this section we reorganize its data
in terms of $\mc{O}_M$-modules using the notion of vertex algebroids
introduced in \cite{GMS2}.

\subsection{} \label{sec.Void}
Given a vertex algebra $V=\bigoplus_{k\geq 0} V_k$, 
the part of its structure involving only $V_0$ and $V_1$ consists of 
the following element and maps
\beq \label{V0.V1}
1\in V_0,\qquad
\d:V_0\rightarrow V_1,\qquad
{}_{(i+j-k-1)}:V_i\times V_j\rightarrow V_k\quad
  \t{for }i,j,k=0\t{ or }1
\eeq
satisfying a number of identities.
Consider ${}_{(-1)}:V_0\times V_0\rightarrow V_0$ and 
${}_{(-1)}:V_0\times V_1\rightarrow V_1$ for example.
The first map makes $A:=V_0$ a commutative ring with unit $1$.
The second map does not make $V_1$ an $A$-module but induces $A$-module
structures on 
\beqq
\Omega:=A_{(-1)}(\d A)\subset V_1 \qquad\t{and}\qquad
\T:=V_1/\Omega.
\eeqq
Choose a splitting $s:\T\rightarrow V_1$ of the projection so as to 
obtain an identification of vector spaces
\beq \label{V1.split}
\Omega\oplus\T\cong V_1,\qquad
(\alpha,X)\mapsto\alpha+s(X).
\eeq
In terms of this identification, ${}_{(-1)}:A\times V_1\rightarrow V_1$ 
takes the form
\beq \label{star.defn}
f_{(-1)}(\alpha,X)=(f\alpha+f\ast X,fX) \quad\t{where}\quad 
f\ast X:=f_{(-1)}s(X)-s(fX).
\eeq
Similarly we can rephrase the other maps in (\ref{V0.V1}) and 
the identities they satisfy in terms of $A$, $\Omega$ and $\T$. 
The resulting data fall into two types.
The first type of data are independent of $s$:
\beq \label{rigid.data}
\begin{array}{l}
\cdot\;
  (A,1)\t{ is a commutative }\bb{C}\t{-algebra with unit and }
  \Omega,\T\t{ are }A\t{-modules}; \\
\cdot\;
  \t{there is an }A\t{-derivation }\d:A\rightarrow\Omega
  \t{ whose image generates }\Omega\t{ as an }A\t{-module}; \\
\cdot\;
  \T\t{ is also a Lie algebra (with Lie bracket }[\;]); \\ 
\cdot\;
  \t{there is an }A\t{-linear Lie algebra homomorphism }
  \T\rightarrow\t{End}\,A\t{ (denoted }X\mapsto X); \\
\cdot\;
  \t{there is a }\bb{C}\t{-linear Lie algebra homomorphism }
  \T\rightarrow\t{End}\,\Omega\t{ (denoted }X\mapsto L_X); \\
\cdot\;
  \t{the }\T\t{-actions on }A\t{ and }\Omega\t{ are }\d\t{-equivariant}; \\
\cdot\;
  \t{the }\T\t{-actions on }A,\,\Omega\t{ and }\T\;(\t{via }X\mapsto[X,-])
  \t{ satisfy the Leibniz rule with respect to} \\
  \qquad A\t{-multiplications}; \\
\cdot\;
  \t{there is an }A\t{-bilinear pairing }
  \langle\;\;\rangle:\Omega\times\T\rightarrow A\t{ satisfying }
  \langle\d f,X\rangle=Xf.
\end{array}
\eeq
The second type of data include three $\bb{C}$-bilinear maps depending
on $s$
\footnote{\label{brackets.defn}
The last two maps are defined by $\{X,Y\}_0=s(X)_{(1)}s(Y)$ and 
$\{X,Y\}_1=s(X)_{(0)}s(Y)-s([X,Y])$.
}
\beq \label{nonrigid.maps}
\ast:A\times\T\rightarrow\Omega,\qquad
\{\;\;\}_0:\T\times\T\rightarrow A,\qquad
\{\;\;\}_1:\T\times\T\rightarrow\Omega,
\eeq
and satisfying the following identities
\beq \label{nonrigid.ids}
\begin{array}{l}
\{X,Y\}_0=\{Y,X\}_0 \qquad\qquad
\d\{X,Y\}_0=\{X,Y\}_1+\{Y,X\}_1 \\
(fg)\ast X-f\ast(gX)-f(g\ast X)=-(Xf)\d g-(Xg)\d f \\
\{X,fY\}_0=f\{X,Y\}_0-\langle f\ast Y,X\rangle-YXf \\
\{X,fY\}_1=f\{X,Y\}_1-L_X(f\ast Y)+(Xf)\ast Y+f\ast[X,Y] \\
X\{Y,Z\}_0-\{[X,Y],Z\}_0-\{Y,[X,Z]\}_0
  =\langle\{X,Y\}_1,Z\rangle+\langle\{X,Z\}_1,Y\rangle \\
L_X\{Y,Z\}_1-L_Y\{X,Z\}_1+L_Z\{X,Y\}_1
  +\{X,[Y,Z]\}_1-\{Y,[X,Z]\}_1-\{[X,Y],Z\}_1 \\
\qquad\qquad
  =\d\langle\{X,Y\}_1,Z\rangle
\end{array}
\eeq
for $f,g\in A$ and $X,Y,Z\in\T$.

In general, any triple $(A,\Omega,\T)$ as in (\ref{rigid.data}) is called
an \emph{extended Lie algebroid}.
The entire collection of data in (\ref{rigid.data})-(\ref{nonrigid.ids}) 
is called a \emph{vertex algebroid}, denoted as 
$(A,\Omega,\T,\ast,\{\;\}_0,\{\;\}_1)$.
According to the above discussion, every vertex algebra (equipped with 
a `splitting') gives rise to a vertex algebroid.

\subsection{} \label{sec.Void.morphism}
Consider a homomorphism of vertex algebras $\Phi:V\rightarrow V'$.
Suppose $(A,\Omega,\T,\ast,\{\;\}_0,\{\;\}_1)$ and
$(A',\Omega',\T',\ast',\{\;\}'_0,\{\;\}'_1)$ are the vertex algebroids
associated to $V$ and $V'$ with respect to some splittings 
$s:\T\rightarrow V_1$ and $s':\T'\rightarrow V'_1$.
Now $\Phi$ induces in the obvious way a map of triples
\beq \label{rigid.morphism}
\varphi:(A,\Omega,\T)\rightarrow(A',\Omega',\T')
\eeq
and, in terms of the identification (\ref{V1.split}), 
$\Phi|_{V_1}:V_1\rightarrow V'_1$ takes the form
\beq \label{Delta.defn}
\Phi(\alpha,X)=\big(\varphi\alpha+\Delta(X),\,\varphi X\big)
\quad\t{where}\quad
\Delta(X):=\Phi s(X)-s'(\varphi X).
\eeq
While (\ref{rigid.morphism}) respects the extended Lie algebroid structures
on the nose, it only respects the rest of the vertex algebroid structures
up to terms linear in $\Delta$, namely
\begin{align}
\varphi f\ast'\varphi X\; &=
  \varphi(f\ast X)+\Delta(fX)-(\varphi f)\Delta(X) \nonumber\\
\{\varphi X,\varphi Y\}'_0 &=
  \varphi\{X,Y\}_0-\langle\Delta(X),\varphi Y\rangle
  -\langle\Delta(Y),\varphi X\rangle  \label{nonrigid.morphism} \\
\{\varphi X,\varphi Y\}'_1 &=
  \varphi\{X,Y\}_1-L_{\varphi X}\Delta(Y)+L_{\varphi Y}\Delta(X)
  -\d\langle\Delta(X),\varphi Y\rangle 
  +\Delta([X,Y]) \nonumber
\end{align}
for $f\in A$ and $X,Y\in\T$. 

The discussion above motivates us to define 
a \emph{morphism of vertex algebroids} as a pair
\beqq
(\varphi,\Delta):
  (A,\Omega,\T,\ast,\{\;\}_0,\{\;\}_1)\rightarrow
  (A',\Omega',\T',\ast',\{\;\}'_0,\{\;\}'_1)
\eeqq
consisting of a map of triples as in (\ref{rigid.morphism}) that respects 
the extended Lie algebroid structures, as well as a map 
$\Delta:\T\rightarrow\Omega'$ that satisfies (\ref{nonrigid.morphism}).
Composition of morphisms is given by
\beq \label{VAoid.mor.comp}
(\varphi',\Delta')\circ(\varphi,\Delta)
  =(\varphi'\varphi,\,\varphi'\Delta+\Delta'\varphi).
\eeq
Together with \S\ref{sec.Void}, this completes the description of 
the category of vertex algebroids $\mf{Void}$ and a forgetful functor 
$G:\mf{V}'\rightarrow\mf{Void}$, where $\mf{V}'$ is the category of 
vertex algebras (equipped with `splittings').
\footnote{
A morphism in $\mf{V}'$ between two objects $V$, $V'$ is simply 
a vertex algebra homomorphism between them.
}

\subsection{} \label{sec.VA.VAoid.functors}
The forgetful functor $G:\mf{V}'\rightarrow\mf{Void}$ has a left adjoint 
$F:\mf{Void}\rightarrow\mf{V}'$.
The composition $GF$ is naturally isomorphic to the identity on 
$\mf{Void}$.
Given a vertex algebroid $(A,\Omega,\T,\ast,\{\;\}_0,\{\;\}_1)$, let us 
briefly describe its \emph{freely generated vertex algebra} 
$V=F(A,\Omega,\T,\ast,\{\;\}_0,\{\;\}_1)$.

The weight-zero and -one components of $V$ are
\beqq
V_0=A,\qquad V_1=\Omega\oplus\T.
\eeqq
Throughout this discussion, we always have $f,g\in A$, $\alpha\in\Omega$ 
and $X,Y\in\T$, regarded as elements of $V$.
For any $u\in V$, denote its vertex operator by $u(z)$.
In particular, $1(z)=\t{id}$ and $(\d f)(z)=\d_z f(z)$.
The vertex operators of $A$ and $\Omega$ commute among themselves and with 
each other.
Normally ordered products with $f(z)$ are given by
\beqq
f(z)g(z)=(fg)(z),\qquad
f(z)\alpha(z)=(f\alpha)(z),\qquad
:f(z)X(z):\;=(fX)(z)+(f\ast X)(z). 
\eeqq
On the other hand, OPEs with $X(z)$ have the following singular parts
\beqq
X(z)f(w) &\sim& \frac{(Xf)(w)}{z-w} \\
X(z)\alpha(w) &\sim& \frac{\langle\alpha,X\rangle(w)}{(z-w)^2}
  +\frac{(L_X\alpha)(w)}{z-w} \\
X(z)Y(w) &\sim& \frac{\{X,Y\}_0(w)}{(z-w)^2}
  +\frac{[X,Y](w)+\{X,Y\}_1(w)}{z-w}
\eeqq
The vertex algebroid axioms guarantee the consistency of the above
relations, i.e.~their combinations do not lead to any nontrivial 
constraints.
Subject to these relations, the fields $f(z)$, $\alpha(z)$, $X(z)$ 
generate $V$ as a vertex algebra.

\subsection{} \label{sec.CDO.VAoid}
Consider the vertex algebra $\CDO_{\t{alg}}(\bb{C}^d)$ defined in 
\S\ref{sec.alg.CDO}. 
Its induced extended Lie algebroid is
\beqq
\big(\{\t{algebraic functions on }\bb{C}^d\},
  \{\t{algebraic }1\t{-forms on }\bb{C}^d\},
  \{\t{algebraic vector fields on }\bb{C}^d\}\big)
\eeqq
equipped with the usual differential on functions, 
Lie bracket on vector fields, 
Lie derivation of functions and $1$-forms by vector fields,
and pairing between $1$-forms and vector fields.
Choose the splitting
\beqq
s:\{\t{algebraic vector fields on }\bb{C}^d\}\rightarrow
  \CDO_{\t{alg}}(\bb{C}^d)_1,\qquad
X\mapsto a_{i,-1}X^i
\eeqq
where $X^i$ are the components of $X$, i.e.~$X=X^i\d_i$.
The three maps in (\ref{nonrigid.maps}) are then given by
\beq \label{CDO.nonrigid.maps}
f\ast X=-X^i\d_i\d_j f\,db_0^j,\quad
\{X,Y\}_0=-(\d_j X^i)(\d_i Y^j),\quad
\{X,Y\}_1=-(\d_k\d_j X^i)(\d_i Y^j)\,db_0^k.
\eeq
This vertex algebroid freely generates $\CDO_{\t{alg}}(\bb{C}^d)$. 
Now if we replace the algebraic functions, $1$-forms and vector fields on
$\bb{C}^d$ by holomorphic ones on an open set $U\subset\bb{C}^d$, the same 
formulae in (\ref{CDO.nonrigid.maps}) define a new vertex algebroid.
Denote its freely generated vertex algebra by $\CDO(U)$.

Consider an isomorphism of vertex algebras 
$\varphi^*_\xi:\CDO(V)\rightarrow\CDO(U)$ as described 
in \S\ref{sec.CDO.iso}.
Its induced isomorphism of vertex algebroids is given by
\beq \label{CDO.VAoid.iso}
(\varphi^*,\Delta_{\varphi,\xi}):
\big(\mc{O}(V),\Omega^1(V),\T(V),\ast,\{\;\}_0,\{\;\}_1\big)\rightarrow
\big(\mc{O}(U),\Omega^1(U),\T(U),\ast,\{\;\}_0,\{\;\}_1\big)
\eeq
where $\varphi^*$ denotes pullback along a biholomorphism 
$\varphi:U\rightarrow V$ and $\Delta_{\varphi,\xi}$ sends $X\in\T(V)$ to
\beq \label{CDO.Delta}
\Delta_{\varphi,\xi}(X)=
  -\d_i(\varphi^*X)^j\,(\theta_\varphi)^i_{\phantom{i}j}
  -\frac{1}{2}\Tr[\theta_\varphi(\varphi^*X)\cdot\theta_\varphi]
  -\frac{1}{2}\,\iota_{\varphi^*X}\xi
\quad\in\Omega^1(U).
\eeq
Recall the notations in \S\ref{sec.CDO.iso}.

\subsection{} \label{sec.presheaf.VAoid}
Given a presheaf of vertex algebras $\mc{V}=\bigoplus_{k\geq 0}\mc{V}_k$, 
the part of its structure involving only $\mc{V}_0$ and $\mc{V}_1$ is 
equivalent to the following set of data:
For each open set $U$, there is a vertex algebroid
\beq \label{presheaf.VAoid.0}
(\mc{A}(U),\Omega(U),\T(U),\ast_U,\{\;\}_{0,U},\{\;\}_{1,U})
\eeq
where a splitting $s_U:\T(U)\rightarrow\mc{V}_1(U)$ has been chosen to
obtain an identification of vector spaces
\beqq
\Omega(U)\oplus\T(U)\cong\mc{V}_1(U),\qquad
(\alpha,X)\mapsto\alpha+s_U(X)
\eeqq
and define $\ast_U,\{\;\}_{0,U},\{\;\}_{1,U}$ e.g.~as in (\ref{star.defn}).
Given $V\subset U$, there is a morphism of vertex algebroids
\beq \label{presheaf.VAoid.1}
(\varphi_{V,U},\Delta_{V,U}):
  (\mc{A}(U),\Omega(U),\T(U),\cdots)\rightarrow
  (\mc{A}(V),\Omega(V),\T(V),\cdots)
\eeq
where $\Delta_{V,U}:\T(U)\rightarrow\Omega(V)$ is defined in terms of 
$s_U$, $s_V$ as in (\ref{Delta.defn}).
Given $W\subset V\subset U$, we have
\beq \label{presheaf.VAoid.2}
\varphi_{W,U}=\varphi_{W,V}\varphi_{V,U},\qquad
\Delta_{W,U}=\varphi_{W,V}\Delta_{V,U}+\Delta_{W,V}\varphi_{V,U}
\eeq
according to the composition law (\ref{VAoid.mor.comp}). 
In general, any collection of data as described in 
(\ref{presheaf.VAoid.0})-(\ref{presheaf.VAoid.2}) is called 
a \emph{presheaf of vertex algebroids}.

The structures of the extended Lie algebroids 
$\big(\mc{A}(U),\Omega(U),\T(U)\big)$ are respected on the nose by 
the restriction maps $\varphi_{V,U}$ and 
hence define morphisms of presheaves.
For example, there are morphisms of presheaves 
(of $\bb{C}$-vector spaces)
\beq \label{presheaf.ext.Lie.ex}
\mc{A}\times\mc{A}\rightarrow\mc{A},\qquad
\mc{A}\times\Omega\rightarrow\Omega,\qquad
\mc{A}\times\T\rightarrow\T
\eeq
making $\mc{A}$ a presheaf of unital commutative $\bb{C}$-algebras and 
$\Omega$, $\T$ both $\mc{A}$-modules.
On the other hand, in view of (\ref{nonrigid.morphism}), the maps 
$\ast_U,\{\;\}_{0,U},\{\;\}_{1,U}$ from various $U$ do not collaborate 
to define morphisms of presheaves;
their `global' meanings are not yet clear.

\subsection{} \label{sec.Cech.conventions} 
Before continuing, we introduce some shorthand notations for \v Cech
complexes.
Given a presheaf of $\bb{C}$-vector spaces $\mc{S}$, let 
$\check C^*(\mc{S})$ be the \v Cech complex of a fixed open cover
$\mf{U}=\{U_\alpha\}_{\alpha\in I}$ valued in $\mc{S}$.

Subscripts like `$0\cdots p$' will stand for a nonempty open set of
the form 
$U_{\alpha_0\cdots\alpha_p}:=U_{\alpha_0}\cap\cdots\cap U_{\alpha_p}$, 
where $\alpha_0,\ldots,\alpha_p\in I$.
For example, given $a\in\check C^p(\mc{S})$, the notation 
`$a_{0\cdots p}$' refers to the value of $a$ on some 
$U_{\alpha_0\cdots\alpha_p}$.
Bilinear maps between \v Cech complexes will always be defined via 
the following Eilenberg-Zilber map 
$\otimesez:\check C^*(\mc{S})\otimes\check C^*(\mc{T})\rightarrow
\check C^*(\mc{S}\otimes\mc{T})$
\beq \label{Eilenberg.Zilber}
(a\otimesez b)_{0\cdots p+q}
  =r^{\mc{S}}_{0\cdots p+q,0\cdots p}(a_{0\cdots p})\otimes
   r^{\mc{T}}_{0\cdots p+q,p\cdots p+q}(b_{p\cdots p+q}),\qquad
a\in\check C^p(\mc{S}),\;
b\in\check C^q(\mc{T})
\eeq
where e.g.~$r^{\mc{S}}_{0\cdots p+q,0\cdots p}$ denotes the restriction map
in $\mc{S}$ associated to an inclusion of the form 
$U_{\alpha_0\cdots\alpha_{p+q}}\subset U_{\alpha_0\cdots\alpha_p}$.
While (\ref{Eilenberg.Zilber}) is strictly associative, it is only graded
symmetric up to homotopy.
This is why the \v Cech complexes studied below carry only
`homotopy versions' of various algebraic structures.

\subsection{} \label{sec.dgVA} 
Continuing with the discussion in \S\ref{sec.presheaf.VAoid}, we consider 
the \v Cech complex $\check C^*(\mc{V})$ valued in a presheaf of vertex 
algebras $\mc{V}$.
Denote the restriction homomorphisms in $\mc{V}$ by $\Phi_{V,U}$. 
The local vacua, i.e.~the vacua of $\mc{V}(U_\alpha)$, constitute a cocycle
in $\check C^0(\mc{V})$. 
The local infinitesimal translation operators form a chain endomorphism
of $\check C^*(\mc{V})$.
The local vertex operators define a bilinear chain map as follows:
Given $a\in\check C^p(\mc{V})$ and $b\in\check C^q(\mc{V})$, let 
$Y(a,z)b\in\check C^{p+q}(\mc{V})[[z,z^{-1}]]$ be given by
\beqq
\big[Y(a,z)b\big]_{0\cdots p+q}=
  Y\big(\Phi_{0\cdots p+q,0\cdots p}(a_{0\cdots p}),z\big)
  \big(\Phi_{0\cdots p+q,p\cdots p+q}(b_{p\cdots p+q})\big).
\eeqq
\footnote{
This is an example of a bilinear map between \v Cech complexes defined 
via the Eilenberg-Zilber map (\ref{Eilenberg.Zilber}).
}
These data satisfy the axioms of a vertex superalgebra except that 
$Y(a,z)b$ may not be in $\check C^*(\mc{V})((z))$ unless $\mf{U}$ 
is finite, and the locality axiom holds only up to homotopy even 
assuming $\mf{U}$ finite.

Let us define a \emph{dg vertex algebra} to be a $4$-tuple 
\beqq
\Big(\;
  V^*=\bigoplus_{k\geq 0}V^*_k,\;
  1\in V^0,\;
  \d:V^*\rightarrow V^*,\;
  Y(\cdot,z)(\cdot):V^*\otimes V^*\rightarrow V^*((z))\;
\Big)
\eeqq
where each $V^*_k$ is a cochain complex, $1$ is a cocycle,
$\d$ and $Y(\cdot,z)(\cdot)$ are degree-preserving chain maps, 
such that the data form a vertex superalgebra.
\footnote{
$V^*=\bigoplus_{k\geq 0}V^*_k$ has two compatible gradings, 
called `degree' and `weight.'
The differential in the chain complex structure has degree one and 
preserves the weight.
The weight is part of the vertex superalgebra structure (\S\ref{sec.VA}).
}
Given a presheaf of vertex algebras $\mc{V}$, the \v Cech complex 
$\check C^*(\mc{V})$ in the case of a finite $\mf{U}$ is an example of 
a homotopy version of a dg vertex algebra.

\subsection{} \label{sec.dgVAoid1}
Given a dg vertex algebra $V^*=\bigoplus_{k\geq 0}V^*_k$, we reorganize 
the part of its structure involving only $V^*_0$ and $V^*_1$.
Regarded simply as a vertex superalgebra, $V^*$ induces 
a vertex superalgebroid
\beqq
(A^*,\Omega^*,\T^*,\ast,\{\;\}_0,\{\;\}_1)
\eeqq
where $A^*:=V^*_0$, $\Omega^*$, $\T^*$ are now cochain complexes and 
a degree-preserving splitting $s:\T^*\rightarrow V^*_1$ has been chosen
to obtain an identification of graded vector spaces
\beq \label{dg.V1.split}
\Omega^*\oplus\T^*\cong V^*_1,\qquad
(\alpha,X)\mapsto\alpha+s(X).
\eeq
The structure maps of the extended Lie superalgebroid 
$(A^*,\Omega^*,\T^*)$ are all chain maps.
Let us call such an object a \emph{dg extended Lie algebroid}.
In terms of (\ref{dg.V1.split}), the differential 
$d:V^*_1\rightarrow V^*_1$ takes the form
\beq \label{dg.Delta.defn}
d(\alpha,X)=(d\alpha+\Delta(X),dX)
\quad\t{where}\quad
\Delta(X):=ds(X)-s(dX).
\eeq
Since $d^2=0$, we have $d\Delta+\Delta d=0$.
\footnote{
In fact, $H^*(\Delta)$ are the connecting maps in the cohomology long 
exact sequence induced by 
$0\rightarrow\Omega^*\rightarrow V^*_1\rightarrow \T^*\rightarrow 0$.
}
On the other hand, ${}_{(-1)}:A^*\otimes V^*_1\rightarrow V^*_1$ takes 
the form
\beq \label{dg.star.defn}
f_{(-1)}(\alpha,X)=(f\alpha+f\ast X,fX)
\quad\t{where}\quad
f\ast X:=f_{(-1)}s(X)-s(fX)
\eeq
just like (\ref{star.defn}).
The fact that ${}_{(-1)}$ is a chain map translates into the equation
(for $f\in A^p$)
\beq \label{dg.star.diff}
d(f\ast X)-(df)\ast X-(-1)^p f\ast(dX)=
-\Delta(fX)+(-1)^p f\Delta(X)
\eeq
i.e.~$\Delta$ is $A^*$-linear up to the homotopy $\ast$.
Similarly $\{\;\}_0,\{\;\}_1$ are also defined using $s$ (see footnote 
\ref{brackets.defn}) and are null homotopies of chain maps constructed 
from $\Delta$, namely (for $X\in\T^p,Y\in\T^q$)
\begin{align}
d\{X,Y\}_0-\{dX,Y\}_0-(-1)^p\{X,dY\}_0 &=
  \langle\Delta(X),Y\rangle+(-1)^{pq}\langle\Delta(Y),X\rangle 
\label{dg.bracket0.diff} \\
d\{X,Y\}_1-\{dX,Y\}_1-(-1)^p\{X,dY\}_1 &=
  (-1)^p L_X\Delta(Y)-(-1)^{(p+1)q}L_Y \Delta(X) \nonumber \\
&\qquad\qquad 
  +\d\langle\Delta(X),Y\rangle
  -\Delta([X,Y])
\label{dg.bracket1.diff}
\end{align}

Motivated by the above discussion, we define 
a \emph{dg vertex algebroid} as a collection of data
\beqq
\Big(A^*,\;\Omega^*,\;\T^*,\;
  \Delta:\T^*\rightarrow\Omega^{*+1},\;
  \ast:A^*\otimes\T^*\rightarrow\Omega^*,\;
  \{\;\}_0:\T^*\otimes\T^*\rightarrow A^*,\;
  \{\;\}_1:\T^*\otimes\T^*\rightarrow\Omega^*\Big)
\eeqq
where $(A^*,\Omega^*,\T^*)$ is a dg extended Lie algebroid, 
$\Delta$ is a degree-one chain map and $\ast$, $\{\;\}_0$, $\{\;\}_1$ 
are degree-preserving maps, such that the data besides $\Delta$ form 
a vertex superalgebroid and the four maps satisfy  
(\ref{dg.star.diff})-(\ref{dg.bracket1.diff}).
The above discussion describes how (resp. a homotopy version of) 
a dg vertex algebra induces (resp. a homotopy version of) 
a dg vertex algebroid.

\subsection{} \label{sec.dgVAoid2}
Consider again a presheaf of vertex algebras $\mc{V}$ and 
its induced presheaf of vertex algebroids as described in 
\S\ref{sec.presheaf.VAoid}.
The local sections $s_U$ provide an identification of graded vector spaces
\beqq
\check C^*(\Omega)\oplus\check C^*(\T)\cong\check C^*(\mc{V}_1),\qquad
(\alpha,X)\mapsto\alpha+s(X)
\eeqq
where $s(X)_{0\cdots p}=s_{0\cdots p}(X_{0\cdots p})$, $p\geq 0$.
Applying the discussion in \S\ref{sec.dgVAoid1} to $V^*=\check C^*(\mc{V})$,
we obtain a homotopy version of a dg vertex algebroid
\beqq
\big(
  \check C^*(\mc{A}),\,\check C^*(\Omega),\,\check C^*(\T),\,
  \Delta,\,\ast,\,\{\;\}_0,\,\{\;\}_1
\big).
\eeqq
Not surprisingly, this homotopy dg vertex algebroid can be expressed  
entirely in terms of the presheaf of vertex algebroids,
i.e.~the data in (\ref{presheaf.VAoid.0})-(\ref{presheaf.VAoid.2}). 
Let us demonstrate this for $\Delta$ and $\ast$.
Denote the restriction maps in $\mc{V}$ by $\Phi_{V,U}$, and those in
$\mc{A}$, $\Omega$, $\T$ by $\varphi_{V,U}$.
Keep in mind the definitions of $\Delta_{V,U}$ and $\ast_U$ according to 
(\ref{Delta.defn}) and (\ref{star.defn}).
Given $X\in\check C^p(\T)$, it follows from (\ref{dg.Delta.defn}) that
\begin{align}
\Delta(X)_{0\cdots p+1}
&=\sum_{i=0}^{p+1}(-1)^i 
  \Phi_{0\cdots p+1,0\cdots\wh{i}\cdots p+1}
  \big(s_{0\cdots\wh{i}\cdots p+1}(X_{0\cdots\wh{i}\cdots p+1})\big)
\nonumber \\
&\qquad\qquad
  -s_{0\cdots p+1}\left(\sum_{i=0}^{p+1}(-1)^i
  \varphi_{0\cdots p+1,0\cdots\wh{i}\cdots p+1}
  (X_{0\cdots\wh{i}\cdots p+1})
  \right)
\nonumber \\
&=\sum_{i=0}^{p+1}(-1)^i
  \Delta_{0\cdots p+1,0\cdots\wh{i}\cdots p+1}
  (X_{0\cdots\wh{i}\cdots p+1}).
\label{Cech.Delta.defn2}
\end{align}
Given $f\in\check C^p(\mc{A})$ and $X\in\check C^q(\T)$, 
it follows from (\ref{dg.star.defn}) that 
\begin{align}
(f\ast X)_{0\cdots p+q}
&=\varphi_{0\cdots p+q,0\cdots p}(f_{0\cdots p})_{(-1)}\,
  \Phi_{0\cdots p+q,p\cdots p+q}
  \big(s_{p\cdots p+q}(X_{p\cdots p+q})\big) 
\nonumber \\
&\qquad\qquad
-s_{0\cdots p+q}\big(
  \varphi_{0\cdots p+q,0\cdots p}(f_{0\cdots p})\,
  \varphi_{0\cdots p+q,p\cdots p+q}(X_{p\cdots p+q})\big) 
\nonumber \\
&=\varphi_{0\cdots p+q,0\cdots p}(f_{0\cdots p})_{(-1)}\,
  \Phi_{0\cdots p+q,p\cdots p+q}
  \big(s_{p\cdots p+q}(X_{p\cdots p+q})\big)
\nonumber \\
&\qquad\qquad
  -\varphi_{0\cdots p+q,0\cdots p}(f_{0\cdots p})_{(-1)}\,
  s_{0\cdots p+q}\big(\varphi_{0\cdots p+q,p\cdots p+q}
  (X_{p\cdots p+q})\big)
\nonumber \\
&\qquad\qquad
  +\varphi_{0\cdots p+q,0\cdots p}(f_{0\cdots p})_{(-1)}\,
  s_{0\cdots p+q}\big(\varphi_{0\cdots p+q,p\cdots p+q}
  (X_{p\cdots p+q})\big)
\nonumber \\
&\qquad\qquad
-s_{0\cdots p+q}\big(
  \varphi_{0\cdots p+q,0\cdots p}(f_{0\cdots p})\,
  \varphi_{0\cdots p+q,p\cdots p+q}(X_{p\cdots p+q})\big)
\nonumber \\
&=\varphi_{0\cdots p+q,0\cdots p}(f_{0\cdots p})\,
  \Delta_{0\cdots p+q,p\cdots p+q}(X_{p\cdots p+q})
\nonumber \\
&\qquad\qquad
  +\varphi_{0\cdots p+q,0\cdots p}(f_{0\cdots p})
  \ast_{0\cdots p+q}
  \varphi_{0\cdots p+q,p\cdots p+q}(X_{p\cdots p+q}).
\label{Cech.star.defn2}
\end{align}
There are similar formulae for $\{\;\}_0$ and $\{\;\}_1$.
\footnote{ \label{Cech.brackets.defn2}
For $X\in\check C^p(\T)$, $Y\in\check C^q(\T)$, we have
\begin{align*}
(\{X,Y\}_0)_{0\cdots p+q}
&=\langle\beta,X'\rangle+\langle\alpha,Y''\rangle
  +\{X',Y''\}_{0,0\cdots p+q} \\
(\{X,Y\}_1)_{0\cdots p+q}
&=L_{X'}\beta-L_{Y''}\alpha+\d\langle\alpha,Y''\rangle
  +\{X',Y''\}_{1,0\cdots p+q}
\end{align*}
where $X'=f^\T_{0\cdots p+q,0\cdots p}(X_{0\cdots p})$,
$Y''=f^\T_{0\cdots p+q,p\cdots p+q}(Y_{p\cdots p+q})$,
$\alpha=\Delta_{0\cdots p+q,0\cdots p}(X_{0\cdots p})$, 
$\beta=\Delta_{0\cdots p+q,p\cdots p+q}(Y_{p\cdots p+q})$.
}
The fact that these formulae define a homotopy dg vertex algebroid can be 
verified using only the structure of the presheaf of vertex algebroids.

\subsection{}
The discussions in the last few subsections can be summarized in 
a commutative diagram:
\beqq
\xymatrix{
  \t{presheaf of vertex algebras}
    \ar[rr]^{\S\ref{sec.dgVA}}\ar[d]_{\S\ref{sec.presheaf.VAoid}} && 
  \t{homotopy dg vertex algebra}
    \ar[d]^{\S\ref{sec.dgVAoid1}} \\
  \t{presheaf of vertex algebroids} 
    \ar[rr]^{\S\ref{sec.dgVAoid2}} &&
  \t{homotopy dg vertex algebroid}
}
\eeqq

\subsection{} \label{sec.CDO.dgVAoid.Cech}
Consider the sheaf of vertex algebras $\CDO_{M,\xi}$ constructed in 
\S\ref{sec.CDO.mfld}.
Let the fixed open cover now be the finite good cover 
$\mf{U}=\{U_1,\ldots,U_N\}$ of $M$ used in \S\ref{sec.CDO.mfld}.
It follows from the goodness of $\mf{U}$ and the filtrations discussed in 
\S\ref{sec.wt.subsheaves} that $\CDO_{M,\xi}$ is acyclic over the open sets 
$U_{\alpha_0\cdots\alpha_p}$. 
Hence $H^*(\CDO_{M,\xi})$ can be computed from the homotopy dg vertex
algebra $\check C^*(\CDO_{M,\xi})$.
On the other hand, since each $\CDO_{M,\xi}(U_{\alpha_0\cdots\alpha_p})$ is
freely generated by a vertex algebroid, all information of
$\check C^*(\CDO_{M,\xi})$ can be recovered from the induced homotopy 
dg vertex algebroid (\S\ref{sec.CDO.VAoid}, \S\ref{sec.dgVAoid2})
\beq \label{dgVAoid.Cech}
\big(
  \check C^*(\mc{O}_M),\,
  \check C^*(\Omega^1_M),\,
  \check C^*(\T_M),\,
  \check\Delta,\,
  \check\ast,\,
  \{\;\check\}_0,\,
  \{\;\check\}_1
\big).
\eeq
This fulfils the goal of reorganizing the data of $\CDO_{M,\xi}$ 
in terms of $\mc{O}_M$-modules.

In order to compute (\ref{dgVAoid.Cech}) explicitly, recall 
the coordinate charts $\varphi_\alpha:U_\alpha\rightarrow\bb{C}^d$ and 
let $W_{\alpha_0\cdots\alpha_p}=
\varphi_{\alpha_0}(U_{\alpha_0\cdots\alpha_p})$.
In subsequent computations, we will
\beq \label{identifications.Cd}
\t{identify }
\left\{\begin{array}{c}
  \CDO_{M,\xi}(U_{\alpha_0\cdots\alpha_p}) \\
  \mc{O}_M(U_{\alpha_0\cdots\alpha_p}) \\
  \Omega^1_M(U_{\alpha_0\cdots\alpha_p}) \\
  \T_M(U_{\alpha_0\cdots\alpha_p}) \\
  \vdots
\end{array}\right\}
\t{ with }
\left\{\begin{array}{c}
  \CDO(W_{\alpha_0\cdots\alpha_p}) \\
  \mc{O}(W_{\alpha_0\cdots\alpha_p}) \\
  \Omega^1(W_{\alpha_0\cdots\alpha_p}) \\
  \T(W_{\alpha_0\cdots\alpha_p}) \\
  \vdots
\end{array}\right\}.
\eeq
For example, the component $X_{0\cdots p}$ of a \v Cech cochain
$X\in\check C^p(\T_M)$ will be understood as a holomorphic vector field
over $W_{\alpha_0\cdots\alpha_p}$.
Consider the restriction homomorphisms in $\CDO_{M,\xi}$
\beqq
\Phi_{0\cdots p,i_0\cdots i_k}=
\begin{cases}
  (\varphi_{i_0 0})^*_{\xi_{i_0 0}}, & i_0>0 \\
  \t{inc}^*_0, & i_0=0
\end{cases}
\eeqq
where $0\leq i_0<\cdots<i_k\leq p$ and `inc' is an inclusion of open 
subsets within $\bb{C}^d$.
According to (\ref{CDO.VAoid.iso}), the induced morphism of vertex 
algebroids is
\beq \label{CDO.presheaf.VAoid.1}
\big(\varphi^*_{0\cdots p,i_0\cdots i_k},\,
\Delta_{0\cdots p,i_0\cdots i_k}\big)=
\begin{cases}
  (\varphi_{i_0 0}^*,\Delta_{i_0 0}), & i_0>0 \\
  (\t{res},0), & i_0=0
\end{cases}
\eeq
where $\Delta_{\beta\alpha}=\Delta_{\varphi_{\beta\alpha},
\xi_{\beta\alpha}}$ is given by (\ref{CDO.Delta}) and $\t{res}=\t{inc}^*$.
In the definition of (\ref{dgVAoid.Cech}), the maps
\beqq
\check\Delta:\check C^*(\T_M)\rightarrow
  \check C^{*+1}(\Omega^1_M),
\qquad
\check\ast:\check C^*(\mc{O}_M)\otimes\check C^*(\T_M)
  \rightarrow\check C^*(\Omega^1_M)
\eeqq
are given by (\ref{Cech.Delta.defn2}) and (\ref{Cech.star.defn2}).
In view of (\ref{CDO.presheaf.VAoid.1}), those formulae now read
\beq 
\label{CDO.Cech.Delta}
&\check\Delta(X)_{0\cdots p+1}=\Delta_{10}(X_{1\cdots p+1}) & \\
\label{CDO.Cech.star}
&(f\,\check\ast\,X)_{0\cdots p+q}= 
  f_{0\cdots p}\ast \varphi_{p0}^*(X_{p\cdots p+q})
  +f_{0\cdots p}\,\Delta_{p0}(X_{p\cdots p+q}) &
\eeq
where $\ast$ is given by (\ref{CDO.nonrigid.maps}).
The other two maps $\{\;\check\}_0,\{\;\check\}_1$ can be computed in 
a similar way.
\footnote{ \label{CDO.Cech.brackets}
It follows from footnote \ref{Cech.brackets.defn2} and 
(\ref{CDO.presheaf.VAoid.1}) that
\begin{align*}
(\{X,Y\check\}_0)_{0\cdots p+q}
&=\{X_{0\cdots p},\varphi_{p0}^*(Y_{p\cdots p+q})\}_0
+\langle\Delta_{p0}(Y_{p\cdots p+q}),X_{0\cdots p}\rangle \\
(\{X,Y\check\}_1)_{0\cdots p+q}
&=\{X_{0\cdots p},\varphi_{p0}^*(Y_{p\cdots p+q})\}_1
+L_{X_{0\cdots p}}[\Delta_{p0}(Y_{p\cdots p+q})]
\end{align*}
where $\{\;\}_0,\{\;\}_1$ (on the right hand sides) are given by 
(\ref{CDO.nonrigid.maps}).
}

\setcounter{equation}{0}
\section{An Infinite Dimensional Dolbeault Complex}

In this section, we construct a sheaf of dg conformal vertex algebras 
which provides a fine resolution of $\CDO_{M,\xi}$ and whose weight-zero 
component is the Dolbeault resolution of $(\CDO_{M,\xi})_0=\mc{O}_M$.

\subsection{} \label{sec.dgVAoid.hom}
Consider a chain homomorphism of dg vertex algebras
$\Phi:(V^*,d)\rightarrow({V'}^*,d')$.
Suppose
\beqq
(A^*,\Omega^*,\T^*,\Delta,\ast,\{\;\}_0,\{\;\}_1),\qquad
({A'}^*,{\Omega'}^*,{\T'}^*,\Delta',\ast',\{\;\}'_0,\{\;\}'_1)
\eeqq
are the dg vertex algebroids associated to $V^*$ and ${V'}^*$
(\S\ref{sec.dgVAoid1}).
Regarded simply as a homomorphism of vertex superalgebras, $\Phi$ induces 
a morphism of vertex superalgebroids (see \S\ref{sec.Void.morphism})
\beqq
(\varphi,h): 
  (A^*,\Omega^*,\T^*,\ast,\{\;\}_0,\{\;\}_1)\rightarrow
  ({A'}^*,{\Omega'}^*,{\T'}^*,\ast',\{\;\}'_0,\{\;\}'_1).
\eeqq
Since $\Phi$ is a chain map, it follows that $\varphi$ is also a chain map
and, by (\ref{Delta.defn}) and (\ref{dg.Delta.defn})
\beq \label{dgDelta.diff}
d'h-hd=\varphi\Delta-\Delta'\varphi
\eeq 
i.e.~$\varphi$ respects $\Delta$ and $\Delta'$ up to the homotopy $h$.

In general, define a \emph{morphism of dg vertex algebroids}
\beqq
(\varphi,h):
  (A^*,\Omega^*,\T^*,\Delta,\ast,\{\;\}_0,\{\;\}_1)\rightarrow
  ({A'}^*,{\Omega'}^*,{\T'}^*,\Delta',\ast',\{\;\}'_0,\{\;\}'_1)
\eeqq
as a morphism between the underlying vertex superalgebroids such that 
$\varphi$ is a chain map and (\ref{dgDelta.diff}) is satisfied.
It is a \emph{(quasi-)isomorphism} if $\varphi$ is a (quasi-)isomorphism.

\subsection{} \label{sec.Cech.Dolbeault}
Our goal is to construct a dg vertex algebroid quasi-isomorphic to 
(\ref{dgVAoid.Cech}) but with the \v Cech complexes replaced by 
Dolbeault complexes.
This will be carried out via \v Cech-Dolbeault complexes.

Given an $\mc{O}_M$-module $\mc{M}$, it admits a fine resolution
\beqq
\xymatrix{
  0 \ar[r] &
  \mc{M} \ar[r] &
  \mc{M}\otimes\E_M \ar[r]^{1\otimes\dbar} &
  \mc{M}\otimes\E^{0,1}_M \ar[r]^{1\otimes\dbar} &
  \mc{M}\otimes\E^{0,2}_M \ar[r] &
  \ldots
}
\eeqq
where $\E^{p,q}_M$ denotes the sheaf of smooth $(p,q)$-forms and tensor
products are taken over $\mc{O}_M$.
Depending on the context, the \v Cech-Dolbeault complex 
$\check C^*(\mc{M}\otimes\E^{0,*}_M)$ will be understood as one of 
the following:
(i) the double complex whose degree-$(p,q)$ term is 
$\check C^p(\mc{M}\otimes\E^{0,q}_M)$ and whose differentials are 
the \v Cech operator $\delta$ and the Dolbeault operator $\dbar$, or 
(ii) the total complex of the double complex just described with 
the differential $D=\delta+(-1)^p\dbar$, where $p$ is the \v Cech degree.
There are quasi-isomorphic embeddings
\beq \label{embeddings.Cech.Dolbeault}
\check C^*(\mc{M})\hookrightarrow
\check C^*(\mc{M}\otimes\E^{0,*}_M)\hookleftarrow
\Gamma(\mc{M}\otimes\E^{0,*}_M)
\eeq
respectively into the degree-$(*,0)$ and -$(0,*)$ terms.

Let $I$ denote a sequence of integers of the form 
$0\leq i_1<\cdots<i_q\leq d$ and 
$d\bar{b}^I=d\bar{b}^{i_1}\wedge\cdots\wedge d\bar{b}^{i_q}$.
\footnote{
From here on, the coordinates of $\bb{C}^d$ are written as $b^i$ 
instead of $b_0^i$.
}
For any smooth $(0,q)$-form $\omega$ on $\bb{C}^d$ (resp.~$U_\alpha$), 
let $\omega_I$ denote its coefficient of $d\bar{b}^I$ 
(resp.~$d\bar{\varphi}_\alpha^I$), i.e.~$\omega=\omega_I d\bar{b}^I$.
We continue to use the shorthand notations introduced in 
\S\ref{sec.Cech.conventions} (with the same open cover used in 
\S\ref{sec.CDO.dgVAoid.Cech}) and the identifications 
(\ref{identifications.Cd}).

\subsection{} \label{sec.dgVAoid.Cech.Dolbeault.1} 
In this subsection, we replace (\ref{dgVAoid.Cech}) with 
a quasi-isomorphic dg vertex algebroid defined using \v Cech-Dolbeault
complexes.
First, we have a sheaf of dg extended Lie algebroids 
$(\mc{O}_M,\Omega^1_M,\T_M)\otimes\E^{0,*}_M$ whose derivation is 
$\d\otimes 1$ and whose bilinear structures are given the bilinear
structures on $(\mc{O}_M,\Omega^1_M,\T_M)$ coupled with the wedge
product on $\E^{0,*}_M$.
For example, the $\T_M\otimes\E^{0,*}_M$-action on $\E^{0,*}_M$ reads
\beqq
Y\eta:=(Y_I\eta_J)\,d\bar{\varphi}_\alpha^I\wedge 
  d\bar{\varphi}_\alpha^J\qquad\t{on }U_\alpha
\eeqq
which is independent of the choice of holomorphic coordinates 
$\varphi_\alpha$.
Applying $\check C^*(-)$ yields a dg extended Lie superalgebroid made up
of \v Cech-Dolbeault complexes.
Then, the last four components of (\ref{dgVAoid.Cech}) extend to maps
between those \v Cech-Dolbeault complexes.
In particular, the two extensions 
\beqq
\check\Delta:\check C^*(\T_M\otimes\E^{0,*}_M)\rightarrow
  \check C^{*+1}(\Omega^1_M\otimes\E^{0,*}_M),\qquad
\check\ast:\check C^*(\E^{0,*}_M)\times\check C^*(\T_M\otimes\E^{0,*}_M)
  \rightarrow\check C^*(\Omega^1_M\otimes\E^{0,*}_M)
\eeqq
are respectively the degree-$(1,0)$ and bidegree-preserving maps given by
the following generalizations of formulae (\ref{CDO.Cech.Delta}) and
(\ref{CDO.Cech.star})
\beq 
&\label{Cech.Dolbeault.Delta1}
(\check\Delta X)_{0\cdots p+1}=
  \Delta_{10}(X_{1\cdots p+1,I})\otimes d\bar{\varphi}_{10}^I & \\
\label{Cech.Dolbeault.star1}
&(\omega\,\check\ast\,X)_{0\cdots p+q}=
  \big[\omega_{0\cdots p,I}\ast\varphi_{p0}^*(X_{p\cdots p+q,J})
  +\omega_{0\cdots p,I}\,\Delta_{p0}(X_{p\cdots p+q,J})\big]\otimes
  d\bar{b}^I\wedge d\bar\varphi_{p0}^J &
\eeq
\footnote{
$\Delta_{\beta\alpha}=\Delta_{\varphi_{\beta\alpha},\xi_{\beta\alpha}}$, 
still defined by (\ref{CDO.Delta}), is now viewed as an operator on 
smooth vector fields.
}
The \v Cech-Dolbeault versions of $\{\;\check\}_0$ and $\{\;\check\}_1$ 
are defined similarly (see footnote \ref{CDO.Cech.brackets}).
It is easy to check that these data constitute a homotopy dg vertex 
algebroid
\beq \label{dgVAoid.Cech.Dolbeault.1}
\big(
  \check C^*(\E^{0,*}_M),\,
  \check C^*(\Omega^1_M\otimes\E^{0,*}_M),\,
  \check C^*(\T_M\otimes\E^{0,*}_M),\,
  \check\Delta,\,
  \check\ast,\,
  \{\;\check\}_0,\,
  \{\;\check\}_1
\big).
\eeq
In particular, $D\check\Delta+\check\Delta D=0$ and equations
(\ref{dg.star.diff})-(\ref{dg.bracket1.diff}) hold with 
$D$, $\check\Delta$, $\check\ast$, $\{\;\check\}_0$, $\{\;\check\}_1$ 
in place of $d$, $\Delta$, $\ast$, $\{\;\}_0$, $\{\;\}_1$.
There is a quasi-isomorphism
\beq \label{qiso1}
(\t{inc},0):
  (\ref{dgVAoid.Cech})\rightarrow(\ref{dgVAoid.Cech.Dolbeault.1})
\eeq
where inc are embeddings of the first type in 
(\ref{embeddings.Cech.Dolbeault}).

\subsection{} \label{sec.dgVAoid.Cech.Dolbeault.2}
Now we would like to construct another homotopy dg vertex algebroid
\beq \label{dgVAoid.Cech.Dolbeault.2}
\big(
  \check C^*(\E^{0,*}_M),\,
  \check C^*(\Omega^1_M\otimes\E^{0,*}_M),\,
  \check C^*(\T_M\otimes\E^{0,*}_M),\,
  \bar\Delta,\,
  \bar\ast,\,
  \{\;\bar\}_0,\,
  \{\;\bar\}_1
\big)
\eeq
such that $\bar\Delta$ has degree $(0,1)$ and there is an isomorphism
\beq \label{qiso2}
(\t{id},h):
  (\ref{dgVAoid.Cech.Dolbeault.1})\rightarrow
  (\ref{dgVAoid.Cech.Dolbeault.2})
\eeq
composed of the identity on the \v Cech-Dolbeault complexes and 
a bidegree-preserving map
\beqq
h:\check C^*(\T_M\otimes\E^{0,*}_M)\rightarrow
  \check C^*(\Omega^1_M\otimes\E^{0,*}_M).
\eeqq
If such an isomorphism exists, it will follow from (\ref{dgDelta.diff}) 
that
\beq \label{Cech.Dolbeault.Deltas}
\check\Delta-\bar\Delta=Dh-hD
\eeq
and from the graded version of (\ref{nonrigid.morphism}) that
(with $\deg X=(p,q)$, $\deg Y=(r,s)$)
\begin{align}
\omega\,\bar\ast\,X\;\; 
& =\;\;\omega\,\check\ast\,X
  +h(\omega\wedge X)-\omega\wedge h(X) 
\nonumber \\
\{X,Y\bar\}_0 
& =\{X,Y\check\}_0-\langle h(X),Y\rangle
  -(-1)^{(p+q)(r+s)}\langle h(Y),X\rangle 
\label{Cech.Dolbeault.nonrigids} \\
\{X,Y\bar\}_1 
& =\{X,Y\check\}_1-L_X h(Y)
  +(-1)^{(p+q)(r+s)}L_Y h(X)
  -\d\langle h(X),Y\rangle
  +h([X,Y])
\nonumber
\end{align}
In fact, we will define $\bar\Delta$, $\bar\ast$, $\{\;\bar\}_0$,  
$\{\;\bar\}_1$ using (\ref{Cech.Dolbeault.Deltas}), 
(\ref{Cech.Dolbeault.nonrigids}) and an appropriate choice of $h$.
It is then a formal consequence of these formulae that 
(\ref{dgVAoid.Cech.Dolbeault.2}) is indeed a homotopy 
dg vertex algebroid. 

\subsection{} \label{sec.h}
The degree-$(1,0)$ and -$(0,1)$ components of equation 
(\ref{Cech.Dolbeault.Deltas}) are
\beq \label{Cech.Dolbeault.Deltas.components}
\check\Delta=\delta h-h\delta,\qquad
\bar\Delta=(-1)^{p+1}(\dbar h-h\dbar)
\eeq
where $p$ is the \v Cech degree.
Assume that $h$ is of the form (with $\deg X=(p,*)$)
\beq \label{Cech.Dolbeault.h}
h(X)_{0\cdots p}=h_0(X_{0\cdots p,I})\otimes d\bar{b}^I
\eeq
for some $\overline{\mc{O}}$-linear first-order differential operators
$h_\alpha$ on $W_\alpha$.
Then by (\ref{Cech.Dolbeault.Delta1}), the first equation in
(\ref{Cech.Dolbeault.Deltas.components}) boils down to
\beq \label{eqn.local.h}
\Delta_{\beta\alpha}=\varphi_{\beta\alpha}^*\circ h_\beta
  -h_\alpha\circ\varphi_{\beta\alpha}^*
\eeq
where $\Delta_{\beta\alpha}=
\Delta_{\varphi_{\beta\alpha},\xi_{\beta\alpha}}$ is defined in
(\ref{CDO.Delta}).
Furthermore, it follows from a computation that any solution to
(\ref{eqn.local.h}) is of the form
\beq \label{h.local}
h_\alpha(Y)=
  \d_i Y^j\,(\Gamma_\alpha)^i_{\phantom{i}j}
  +\frac{1}{2}\Tr[\Gamma_\alpha(Y)\cdot\Gamma_\alpha]
  -\frac{1}{2}\iota_Y B_\alpha\qquad
\t{in }W_\alpha
\eeq
where $\Gamma_\alpha=[(\Gamma_\alpha)^i_{\phantom{i}j}]$ is a smooth,
matrix-valued $(1,0)$-form on $W_\alpha$ such that
\beq \label{eqn.Gamma}
g_{\beta\alpha}^{-1}\cdot\varphi_{\beta\alpha}^*\Gamma_\beta\cdot
  g_{\beta\alpha}
-\Gamma_\alpha
= -\theta_{\beta\alpha}\qquad
\t{in }W_{\alpha\beta}
\eeq
and $B_\alpha$ is a smooth $(2,0)$-form on $W_\alpha$ such that
\beq \label{eqn.B}
\varphi_{\beta\alpha}^*B_\beta
-B_\alpha
= \xi_{\beta\alpha}
  +\Tr[\theta_{\beta\alpha}\wedge\Gamma_\alpha]\qquad
\t{in }W_{\alpha\beta}.
\eeq
The existence of such $\{\Gamma_\alpha\}$ and $\{B_\alpha\}$ is addressed
below.

\subsection{} \label{sec.connection}
A collection $\{\Gamma_\alpha\in\E^1(W_\alpha)\otimes\mf{gl}(d)\}$ 
satisfying (\ref{eqn.Gamma}) is precisely a set of local connection
$1$-forms on $TM$ that are compatible with the coordinate transformations 
$\varphi_{\beta\alpha}$.
Also, if $\{\Gamma_\alpha\}$ satisfy (\ref{eqn.Gamma}), 
so do their $(1,0)$-parts.
This proves the existence of $\{\Gamma_\alpha\}$ with 
the desired properties.

Given a hermitian metric on $M$, there is a unique unitary connection 
on $TM$ whose connection $1$-forms are of pure type $(1,0)$.~\cite{Wells}
This is sometimes called the Chern connection.
From here on, $\nabla$ will always denote the Chern connection with respect
to a fixed hermitian metric, and $\{\Gamma_\alpha\}$ its connection
$(1,0)$-forms.
The curvature form $R$ associated to $\nabla$ is of pure type $(1,1)$.

\subsection{} \label{sec.B}
Consider again the computations in the \v Cech-de Rham complex shown in 
figure \ref{fig.CdR.ch2}.
By (\ref{glue.assoc}), the degree-$(1,2)$ cochain in the diagram is 
$\delta$-closed, hence also $\delta$-exact as each row is acyclic.
In other words, there exist $\{B_\alpha\in\E^2(W_\alpha)\}$ that satisfy 
(\ref{eqn.B}).
Also, if $\{B_\alpha\}$ satisfy (\ref{eqn.B}), so do their $(2,0)$-parts.
This proves the existence of $\{B_\alpha\}$ with the desired properties.

Observe that the cochain $\{dB_\alpha-CS(\Gamma_\alpha)\}$ is 
$\delta$-closed.
Let $H$ be the global $3$-form such that
\beqq
H|_{U_\alpha}=\varphi_\alpha^*(dB_\alpha-CS(\Gamma_\alpha)).
\eeqq
By construction, $H$ has only $(3,0)$- and $(2,1)$-parts, and 
$dH=-\Tr(R\wedge R)$.
In fact, figure \ref{fig.CdR.ch2} shows that the existence of the following
are equivalent:
\footnote{ \label{O2cl.resoln}
This is easy to show directly.
More conceptually, if $H^*(\Omega^{2,\t{cl}}_M)$ is computed using 
the fine resolution
\beqq
\xymatrix{
  0 \ar[r] &
  \Omega^{2,\t{cl}}_M \ar[r] &
  \E^{2,0}_M \ar[r]^{d\phantom{dddd}}  &
  \E^{3,0}_M\oplus\E^{2,1}_M \ar[r]^{d\phantom{dddd}} &
  \E^{4,0}_M\oplus\E^{3,1}_M\oplus\E^{2,2}_M \ar[r]^{\phantom{dddddd}d} &
  \ldots
}
\eeqq
figure \ref{fig.CdR.ch2} shows that the first obstruction class studied 
in \S\ref{sec.obstructions} is represented by the global $(2,2)$-form
$-\Tr(R\wedge R)$.
}
\beqq
\begin{array}{l}
\cdot\;
  \{\xi_{\beta\alpha}\in\E^{2,0}(W_{\alpha\beta})\}
  \t{ satisfying }d\xi_{\beta\alpha}=WZ_{\beta\alpha}
  \t{ and }(\ref{glue.assoc}) 
  \vspace{0.05in}\\
\cdot\;
  \{B_\alpha\in\E^{2,0}(W_\alpha)\}\t{ such that }
  \{dB_\alpha-CS(\Gamma_\alpha)\}\t{ glue into a global }3\t{-form}
  \vspace{0.05in}\\
\cdot\;
  H\in\E^{3,0}(M)\oplus\E^{2,1}(M)\t{ satisfying }dH=-\Tr(R\wedge R)
\end{array}
\eeqq

\subsection{} \label{sec.dgVAoid.Dolbeault}
Now that we have solved the first equation in 
(\ref{Cech.Dolbeault.Deltas.components}), we define the homotopy dg vertex
algebroid (\ref{dgVAoid.Cech.Dolbeault.2}) using 
(\ref{Cech.Dolbeault.nonrigids}) and the second equation in
(\ref{Cech.Dolbeault.Deltas.components}).
Since $\bar\Delta$ has degree $(0,1)$, the equation 
$D\bar\Delta+\bar\Delta D=0$ implies that $\bar\Delta$ anticommutes with 
$\delta$ and thus is induced by a morphism of sheaves
\beqq
\bar\Delta:\T_M\otimes\E_M^{0,*}\rightarrow\Omega^1_M\otimes\E_M^{0,*+1}.
\eeqq
Axioms (\ref{dg.star.diff})-(\ref{dg.bracket1.diff}) for 
(\ref{Cech.Dolbeault.Deltas.components}) similarly imply that 
$\bar\ast,\{\;\bar\}_0,\{\;\bar\}_1$ are induced by morphisms of 
sheaves
\begin{align*}
\bar\ast\;\;&:
  \E^{0,*}_M\times
  \big(\T_M\otimes\E^{0,*}_M\big)\rightarrow
  \Omega^1_M\otimes\E^{0,*}_M \\
\{\;\bar\}_0&:
  \big(\T_M\otimes\E^{0,*}_M\big)\times
  \big(\T_M\otimes\E^{0,*}_M\big)\rightarrow
  \E^{0,*}_M \\
\{\;\bar\}_1&:
  \big(\T_M\otimes\E^{0,*}_M\big)\times
  \big(\T_M\otimes\E^{0,*}_M\big)\rightarrow
  \Omega^1_M\otimes\E^{0,*}_M
\end{align*}
These operators define a sheaf of dg vertex algebroids
\beq \label{sheaf.dgVAoid.Dolbeault}
\big(
  \E^{0,*}_M,\,\Omega^1_M\otimes\E^{0,*}_M,\,\T_M\otimes\E^{0,*}_M,\,
  \bar\Delta,\,\bar\ast,\,\{\;\bar\}_0,\,\{\;\bar\}_1
\big)
\eeq
whose structure is strict, i.e.~not merely up to homotopy.
\footnote{ 
This is expected because e.g.~at weight zero, while the \v Cech complexes 
$\check C^*(\mc{O}_M)$, $\check C^*(\E^{0,*}_M)$ are graded commutative 
only up to homotopy, the Dolbeault complex $\E^{0,*}_M$ is strictly 
graded commutative.
}
There is a quasi-isomorphism
\beq \label{qiso3}
(\t{inc},0):
  \Gamma(\ref{sheaf.dgVAoid.Dolbeault})\rightarrow
  (\ref{dgVAoid.Cech.Dolbeault.2})
\eeq
where inc are embeddings of the second type in
(\ref{embeddings.Cech.Dolbeault}).

Let us obtain formulae for $\bar\Delta$, $\bar\ast$, $\{\;\bar\}_0$,
$\{\;\bar\}_1$.
First restrict to degree $(0,0)$.
By (\ref{Cech.Dolbeault.Deltas.components}) and (\ref{Cech.Dolbeault.h})
\beq \label{Dolbeault0.Delta.local}
\bar\Delta(X)_\alpha
  =[-\bar\d_i h_\alpha(X_\alpha)+h_\alpha(\bar\d_i X_\alpha)]
  \otimes d\bar b^i.
\eeq
On the other hand, by (\ref{Cech.Dolbeault.star1}), 
(\ref{Cech.Dolbeault.nonrigids}), (\ref{Cech.Dolbeault.h}) as well as 
footnote \ref{CDO.Cech.brackets}, we have
\begin{align}
(f\,\bar\ast\,X)_\alpha\;\;&=\; 
  f_\alpha\ast X_\alpha
  +h_\alpha(f_\alpha X_\alpha)
  -f_\alpha h_\alpha(X_\alpha)
\nonumber \\
\big(\{X,Y\bar\}_0\big)_\alpha &= 
  \{X_\alpha,Y_\alpha\}_0
  -\langle h_\alpha(X_\alpha),Y_\alpha\rangle
  -\langle h_\alpha(Y_\alpha),X_\alpha\rangle
\label{Dolbeault0.nonrigids.local} \\
\big(\{X,Y\bar\}_1\big)_\alpha &= 
  \{X_\alpha,Y_\alpha\}_1
  -L_{X_\alpha}h_\alpha(Y_\alpha)
  +L_{Y_\alpha}h_\alpha(X_\alpha)
  -\d\langle h_\alpha(X_\alpha),Y_\alpha\rangle
  +h_\alpha([X_\alpha,Y_\alpha])
\nonumber
\end{align}
Applying (\ref{h.local}) and (\ref{CDO.nonrigid.maps}) yields invariant
expressions for 
(\ref{Dolbeault0.Delta.local})-(\ref{Dolbeault0.nonrigids.local}).
To write them down, pick a local frame $e_1,\ldots,e_d$ of $TM$,
let $e^1,\ldots,e^d$ be its dual frame of $T^*M$, and associate to every
vector field $X$ a section of $\t{End}\,TM$ as follows
\beqq
\wt\nabla X:=\nabla X+T(X,-).
\eeqq
Also, repeated indices will be summed over.
The invariant expression for (\ref{Dolbeault0.Delta.local}) is given by
\beq \label{Dolbeault0.Delta}
\bar\Delta(X)
=\left[\Tr(\wt\nabla X\circ R_{e_i,\bar e_j})
  +\frac{1}{2}H(X,e_i,\bar e_j)\right]\,e^i\otimes\bar e^j
\eeq
where $R$ is the curvature $(1,1)$-form and $H$ is the $3$-form constructed
in \S\ref{sec.B}.
The invariant expressions for (\ref{Dolbeault0.nonrigids.local}) are 
given by
\begin{align} 
f\,\bar\ast\,X\;\;\;&=
  -(\nabla_{e_i}\d f)(X)\,e^i
\nonumber \\
\{X,Y\bar\}_0&=
  -\Tr(\wt\nabla X\circ\wt\nabla Y) 
\label{Dolbeault0.nonrigids} \\
\{X,Y\bar\}_1&=
  \left[-\Tr\big(\nabla_{e_i}(\wt\nabla X)\circ\wt\nabla Y\big)
  +\frac{1}{2}H(X,Y,e_i)\right]e^i
\nonumber
\end{align} 
Now consider higher Dolbeault degrees.
For a smooth function $f$, smooth vector fields $X,Y$, and 
antiholomorphic $(0,*)$-forms $\omega,\eta$, we have
\beq \label{Dolbeault.all}
\begin{array}{ll}
\bar\Delta(X\otimes\omega)
  =\bar\Delta(X)\wedge\omega &
\{X\otimes\omega,Y\otimes\eta\bar\}_0
  =\{X,Y\bar\}_0\,\omega\wedge\eta
\vspace{0.05in} \\
(f\omega)\,\bar\ast\,(X\otimes\eta)
  =(f\,\bar\ast\,X)\otimes(\omega\wedge\eta)\qquad &
\{X\otimes\omega,Y\otimes\eta\bar\}_1
  =\{X,Y\bar\}_1\otimes(\omega\wedge\eta) 
\end{array}
\eeq
which are well-defined because 
(\ref{Dolbeault0.Delta})-(\ref{Dolbeault0.nonrigids}) are 
$\overline{\mc{O}}_M$-linear.

According to (\ref{Dolbeault0.nonrigids.local}), there are morphisms of 
sheaves of vertex algebroids
\beq \label{dgVAoid.mor.CDO.Dolbeault}
\big(\varphi_\alpha^*,\varphi_\alpha^*h_\alpha\big):
\big(\mc{O}_{W_\alpha},\Omega^1_{W_\alpha},\T_{W_\alpha},
  \ast,\{\;\}_0,\{\;\}_1\big)\rightarrow
\big(\E_M,\,\Omega^1_M\otimes\E_M,\,\T_M\otimes\E_M,\,
  \bar\ast,\,\{\;\bar\}_0,\,\{\;\bar\}_1\,\big)
\eeq
Furthermore, it follows from (\ref{eqn.local.h}) that there are 
commutative diagrams
\beq \label{dgVAoid.mor.CDO.CDO.Dolbeault}
\qquad
\xymatrix{
  \big(\mc{O}_{W_{\alpha\beta}},
    \Omega^1_{W_{\alpha\beta}},
    \T_{W_{\alpha\beta}},
    \ast,\{\;\}_0,\{\;\}_1\big)
  \ar[rd]^{\phantom{aaa}(\varphi_\alpha^*,\,\varphi_\alpha^*h_\alpha)} 
  & \\
  &\qquad
  \big(\E_M,\,\Omega^1_M\otimes\E_M,\,\T_M\otimes\E_M,\,
    \bar\ast,\,\{\;\bar\}_0,\,\{\;\bar\}_1\big) \\
  \big(\mc{O}_{W_{\beta\alpha}},
    \Omega^1_{W_{\beta\alpha}},
    \T_{W_{\beta\alpha}},
    \ast,\{\;\}_0,\{\;\}_1\big)
  \ar[ru]_{\phantom{aaa}(\varphi_\beta^*,\,\varphi_\beta^*h_\beta)} 
  \ar[uu]^{(\varphi_{\beta\alpha}^*,\,\Delta_{\beta\alpha})} &&
}
\eeq

\subsection{} \label{sec.dgVA.Dolbeault}
Consider the sheaf of dg vertex algebras freely generated by 
(\ref{sheaf.dgVAoid.Dolbeault})
\beq \label{sheaf.dgVA.Dolbeault}
\big(\CDC_{M,H},\dbar^{\t{ch}}\big):=
F\big(
  \E^{0,*}_M,\,\Omega^1_M\otimes\E^{0,*}_M,\,\T_M\otimes\E^{0,*}_M,\,
  \bar\Delta,\,\bar\ast,\,\{\;\bar\}_0,\,\{\;\bar\}_1
\big).
\eeq
For a description of its vertex superalgebra structure, 
see \S\ref{sec.VA.VAoid.functors}.
The differential at weight zero is
\beqq
\dbar^{\t{ch}}_0=\dbar\;\;\t{on}\;\E^{0,*}_M.
\eeqq
On the other hand, by (\ref{dg.Delta.defn}) the differential at weight one
is a `deformed' Dolbeault operator
\beq \label{dgVA.Dolbeault.1}
\dbar^{\t{ch}}_1=\dbardef\;\t{on}\;\;
  (\Omega^1_M\oplus\T_M)\otimes\E^{0,*}_M,
\qquad
(\alpha,X)\mapsto\big(\dbar\alpha+\bar\Delta(X),\,\dbar X\big).
\eeq
Notice that (\ref{sheaf.dgVA.Dolbeault}) depends on $H$ via the definition 
of $\bar\Delta$ and $\{\;\bar\}_1$ in 
(\ref{Dolbeault0.Delta})-(\ref{Dolbeault.all}).
Recall the sequence of quasi-isomorphisms (\ref{qiso1}), (\ref{qiso2}) 
and (\ref{qiso3}).
Since $\check C^*(\CDO_{M,\xi})$ and $\Gamma(\CDC_{M,H})$ are freely 
generated by quasi-isomorphic dg vertex algebroids, we may expect that 
they are quasi-isomorphic dg vertex algebras and hence both compute 
$H^*(\CDO_{M,\xi})$.
This will be argued more carefully below.

Applying the free functor $F$ to (\ref{dgVAoid.mor.CDO.Dolbeault}) and 
(\ref{dgVAoid.mor.CDO.CDO.Dolbeault}) yields morphisms of sheaves 
of vertex algebras 
\beq \label{inc.CDO.CDC}
i_\alpha:=F(\varphi_\alpha^*,\varphi_\alpha^*h_\alpha):
  \CDO|_{W_\alpha}\rightarrow\E^{\t{ch},0}_{M,H}
\eeq
and commutative diagrams of the form (recall the last paragraph of 
\S\ref{sec.CDO.VAoid})
\beqq
\xymatrix{
  \CDO|_{W_{\alpha\beta}}\ar[rrd]^{i_\alpha} && \\
  && \E^{\t{ch},0}_{M,H}\;. \\
  \CDO|_{W_{\beta\alpha}}
    \ar[uu]^{(\varphi_{\beta\alpha})^*_{\xi_{\beta\alpha}}}
    \ar[rru]_{i_\beta} &&
}
\eeqq
This is equivalent to a morphism of sheaves of vertex algebras 
$i:\CDO_{M,\xi}\rightarrow\E^{\t{ch},0}_{M,H}$ (\S\ref{sec.CDO.mfld}).
Composition with the inclusion 
$\E^{\t{ch},0}_{M,H}\hookrightarrow\CDC_{M,H}$ then defines 
a fine resolution
\beq \label{chiral.Dolbeault.resoln}
\xymatrix{
  0\ar[r] & \CDO_{M,\xi}\ar[r]^{i\phantom{aaa}} & 
  (\CDC_{M,H},\dbar^{\t{ch}}).
}
\eeq
Indeed, the weight-zero component, namely the ordinary Dolbeault 
resolution, is exact. 
The weight-one component, which is described in (\ref{dgVA.Dolbeault.1}), 
fits into the commutative diagram in figure \ref{fig.resoln.wt1.ses}
in which exactness of the top row, bottom row and all columns implies 
that of the middle row.
Applying this argument repeatedly to filtrations like 
(\ref{wt2.filtration}) proves exactness at all higher weights.
The resolution (\ref{chiral.Dolbeault.resoln}) leads to the following 
isomorphism of vertex superalgebras
\beq \label{dgVA.Dolbeault.chlgy}
H^*\big(\Gamma(\CDC_{M,H}),\dbar^{\t{ch}}\big)\cong
H^*(\CDO_{M,\xi}).
\eeq

\begin{figure}[t] 
{\scriptsize
\beqq
\xymatrix{
  & 0\ar[d] && 0\ar[d] && 0\ar[d] && \\
  0\ar[r] & \Omega^1_M\ar[rr]\ar[d] && 
    \Omega^1_M\otimes\E_M\ar[rr]^{\bar\d}\ar[d] &&
    \Omega^1_M\otimes\E^{0,1}_M\ar[rr]^{\bar\d}\ar[d] && \cdots \\
  0\ar[r] & (\CDO_{M,\xi})_1\ar[rr]^i\ar[d] &&
    (\Omega^1_M\oplus\T_M)\otimes\E_M\ar[rr]^{\dbardef}\ar[d] &&
    (\Omega^1_M\oplus\T_M)\otimes\E^{0,1}_M
    \ar[rr]^{\phantom{aaa}\dbardef}\ar[d] && 
    \cdots \\
  0\ar[r] & \T_M\ar[rr]\ar[d] && 
    \T_M\otimes\E_M\ar[rr]^{\bar\d}\ar[d] &&
    \T_M\otimes\E^{0,1}_M\ar[rr]^{\bar\d}\ar[d] && \cdots \\
  & 0 && 0 && 0 &&
}
\eeqq
}
\caption{
  The two-step filtration of the resolution 
  (\ref{chiral.Dolbeault.resoln}) at weight one.
  \label{fig.resoln.wt1.ses}}
\end{figure}

\subsection{} \label{sec.CDC.iso.classes}
Consider an isomorphism of sheaves of dg vertex algebras between 
$\CDC_{M,H}$ and $\CDC_{M,H'}$ over the identity on $M$.
This is equivalent to an isomorphism between the associated sheaves of 
dg vertex algebroids
\beqq
&&
(\t{id},\beta):
\big(
  \E^{0,*}_M,\,
  \Omega^1_M\otimes\E^{0,*}_M,\,
  \T_M\otimes\E^{0,*}_M,\,
  \bar\Delta,\,
  \bar\ast,\,
  \{\;\bar\}_0,\,
  \{\;\bar\}_1
\big) \\
&&
\phantom{(\t{id},\t{id},\t{id},\rho):}\qquad
\longrightarrow
\big(
  \E^{0,*}_M,\,
  \Omega^1_M\otimes\E^{0,*}_M,\,
  \T_M\otimes\E^{0,*}_M,\,
  \bar\Delta',\,
  \bar\ast,\,
  \{\;\bar\}_0,\,
  \{\;\bar\}'_1
\big)
\eeqq
where $\beta:\T_M\otimes\E^{0,*}_M\rightarrow\Omega^1_M\otimes\E^{0,*}_M$
is a degree-preserving operator.
According to (\ref{dgDelta.diff}) and (\ref{nonrigid.morphism}), 
$\beta$ has to satisfy precisely the following conditions 
\beq
& \dbar[\beta(X)]-\beta(\dbar X)=\bar\Delta(X)-\bar\Delta'(X) & 
\label{iso.Delta}  \\
& \beta(\omega\wedge X)-\omega\wedge\beta(X)=0 & 
\label{iso.star} \\
& \langle\beta(X),Y\rangle+(-1)^{pq}\langle\beta(Y),X\rangle=0 & 
\label{iso.bracket0} \\
& L_X\beta(Y)-(-1)^{pq}L_Y\beta(X)+\d\langle\beta(X),Y\rangle-\beta([X,Y])
  =\{X,Y\bar\}_1-\{X,Y\bar\}'_1 &
\label{iso.bracket1}
\eeq
where $\omega,X,Y$ are respectively sections of $\E^{0,*}_M$, 
$\T_M\otimes\E^{0,p}_M$ and $\T_M\otimes\E^{0,q}_M$.
By (\ref{iso.star}), $\beta$ is determined by its component in degree zero.
Then (\ref{iso.star}) and (\ref{iso.bracket0}) together imply that
\beqq
\tilde\beta(X,Y):=\langle\beta(X),Y\rangle\quad 
\t{for vector fields }X,Y
\eeqq
is a smooth $(2,0)$-form.
Using formulae (\ref{Dolbeault0.Delta})-(\ref{Dolbeault0.nonrigids}),
conditions (\ref{iso.Delta}) and (\ref{iso.bracket1}) can be rewritten as
\beqq
\bar\d\tilde\beta=\frac{1}{2}(H-H')^{2,1},\qquad
\d\tilde\beta=\frac{1}{2}(H-H')^{3,0}.
\eeqq
Hence $\CDC_{M,H},\CDC_{M,H'}$ are isomorphic over the identity
if and only if $H,H'$ differ by the de Rham deriviative of 
a smooth $(2,0)$-form.
\footnote{
Therefore by the resolution in footnote \ref{O2cl.resoln}, the isomorphism 
classes of $\{\CDC_{M,H}\}$ form an affine space modeled on 
$H^1(\Omega^{2,\t{cl}}_M)$. 
See also \S\ref{sec.CDO.iso.classes}.
}

\subsection{} \label{sec.CDC.conformal}
Now we address the conformal structure of $\CDC_{M,H}$.
First let us rephrase condition (\ref{glue.conformal}) in global terms.
Recall \S\ref{sec.connection}.
According to figure \ref{fig.CdR.ch1}, (\ref{glue.conformal}) implies 
that the cochain $\{\Tr\Gamma_\alpha\}$ is $\delta$-closed.
Let $A$ be the global $(1,0)$-form such that
\beqq
A|_{U_\alpha}=\varphi_\alpha^*(\Tr\Gamma_\alpha).
\eeqq
Then $dA=\Tr R$.
In fact, the existence of the following are equivalent:
\footnote{ 
If $H^*(\Omega^{1,\t{cl}}_M)$ is computed using the fine resolution
\beqq
\xymatrix{
  0 \ar[r] &
  \Omega^{1,\t{cl}}_M \ar[r] &
  \E^{1,0}_M \ar[r]^{d\phantom{dddd}}  &
  \E^{2,0}_M\oplus\E^{1,1}_M \ar[r]^{\phantom{dddd}d} &
  \ldots
}
\eeqq
figure \ref{fig.CdR.ch1} shows that the second obstruction class studied 
in \S\ref{sec.obstructions} is represented by the global $(1,1)$-form
$\Tr R$.
}
\beqq
\begin{array}{l}
\cdot\;
  \t{coordinate charts }\varphi_\alpha\t{ satisfying }(\ref{glue.conformal}) 
  \vspace{0.05in}\\
\cdot\;
  \t{coordinate charts }\varphi_\alpha\t{ such that }
  \{\Tr\Gamma_\alpha\}\t{ glue into a global form}
  \vspace{0.05in}\\
\cdot\;
  A\in\E^{1,0}(M)\t{ satisfying }dA=\Tr R
\end{array}
\eeqq

For each $\alpha$, denote by $\nu_\alpha$ the local conformal element 
(\ref{conformal}) of $\CDO_{M,\xi}(U_\alpha)=\CDO(W_\alpha)$.
By a computation using (\ref{inc.CDO.CDC}) and (\ref{h.local}), the image 
of $\nu_\alpha$ in $\E^{\t{ch},0}_{M,H}(U_\alpha)$ equals 
\beq \label{i.conformal}
i_\alpha(\nu_\alpha)
=\sum_{i=1}^d
  \Big(\frac{\d}{\d\varphi_\alpha^i}\Big)_{-1}(d\varphi_\alpha^i)
+\frac{1}{2}\Tr\big[(\varphi_\alpha^*\Gamma_\alpha)_{-1}
  (\varphi_\alpha^*\Gamma_\alpha)\big].
\eeq
It can be checked that (\ref{i.conformal}) is a conformal element of 
the larger vertex superalgebra $\CDC_{M,H}(U_\alpha)$ with the same 
central charge $2d$.
By assumption, as $\alpha$ varies, the local sections (\ref{i.conformal}) 
glue into a global one.
Denote this global section as well as any of its restrictions by $\nu$.

Consider the Virasoro field $L(z)$ associated to $\nu$.
Let $f$ be a smooth function and $X$ a smooth vector field on $M$.
OPEs with $L(z)$ have the following singular parts
\beqq
L(z)f(w) &\sim& \frac{\d_w f(w)}{z-w} \\
L(z)X(w) &\sim& 
  \frac{[\Tr\wt\nabla X-A(X)](w)}{(z-w)^3}
  +\frac{X(w)}{(z-w)^2}
  +\frac{\d_w X(w)}{z-w}
\eeqq
Since $\nu$ is in the image of the morphism $i$ in 
(\ref{chiral.Dolbeault.resoln})
\beqq
\dbar^{\t{ch}}(\nu)=0
\quad\Rightarrow\quad
[\dbar^{\t{ch}},L(z)]=0.
\eeqq
Therefore $\nu$ induces a conformal structure on the cohomology of 
$\big(\Gamma(\CDC_{M,H}),\dbar^{\t{ch}}\big)$, making 
(\ref{dgVA.Dolbeault.chlgy}) an isomorphism of 
conformal vertex superalgebras.
By (\ref{CDO.mfld.Euler}), we have
\beq \label{CDC.global.chlgy}
\t{char}\,H^*\big(\Gamma(\CDC_{M,H}),\dbar^{\t{ch}}\big)
=\frac{W(M)}{\eta(q)^{2d}}
\eeq
which provides a new geometric interpretation of the Witten genus.

\subsection{} \label{sec.CDC.Kahler}
Let us finish with a comment on the geometric conditions stated
in \S\ref{sec.B} and \S\ref{sec.CDC.conformal}:
\beqq
\begin{array}{ll}
\t{(i)}
&\exists\,H\in\E^{3,0}(M)\oplus\E^{2,1}(M)\t{ such that }
  dH=-\Tr(R\wedge R) 
\vspace{0.02in} \\
\t{(ii)}
&\exists\,A\in\E^{1,0}(M)\t{ such that }dA=\Tr R
\end{array}
\eeqq
Recall that (i) is required for the construction of either $\CDO_{M,\xi}$
or $\CDC_{M,H}$, and (ii) is required for a global conformal structure.
In the case $M$ is K\"ahler, they are equivalent to the following
topological conditions
\beqq
\begin{array}{ll}
\t{(i)}'
&\exists\,H\in\E^3(M)\t{ such that }dH=-\Tr(R\wedge R)\phantom{WWWWWW}
\vspace{0.02in} \\
\t{(ii)}'
&\exists\,A\in\E^1(M)\t{ such that }dA=\Tr R
\end{array}
\eeqq
For example, assume (ii)$'$. 
Since $R$ has no $(0,2)$-part, $\dbar A^{0,1}=0$.
Suppose $\omega$ is the Dolbeault harmonic representative of 
$[A^{0,1}]\in H^{0,1}(M)$.
Hence $A^{0,1}-\omega=\dbar f$ for some $f\in\E(M)$.
Let $A'=A-\omega-df$.
By construction $A'$ has no $(0,1)$-part.
Since $\omega$ is also de Rham harmonic, $dA'=dA=\Tr R$.
This proves (ii).
The proof of (i)$'\Rightarrow$ (i) is similar.
By Chern-Weil theory, (i)$'$ and (ii)$'$ are equivalent to the vanishing 
of $ch_2(M)$ and $ch_1(M)$ respectively.

{\footnotesize

\vspace{0.1in}
\rule{2.5in}{.1pt}
\vspace{0.1in}

\noindent
{\sc Max-Planck-Institut f\"ur Mathematik, 53111 Bonn, Germany} \\
{\it Email address:}
{\sf pokman@mpim-bonn.mpg.de}

}

\end{document}